\documentclass[10pt, letterpaper]{article}
\usepackage[utf8]{inputenc}

\usepackage{microtype}
\usepackage{hyperref}
\usepackage{authblk}
\usepackage[dvipsnames]{xcolor}
\usepackage[caption=false,font=normalsize,labelfont=sf,textfont=sf]{subfig}
\usepackage{graphicx}
\usepackage{amsmath,amssymb,amsfonts,amsthm} 
\usepackage[mathscr]{eucal}

\newtheorem{proposition}{Proposition}[section]

\newtheorem{ass}{Assumption}
\newtheorem{example}{Example}
\newcommand{\R}{\mathbb{R}}
\title{Coexistence of hidden attractors and self-excited attractors through breaking heteroclinic-like orbits of switched systems}

\author{R.J.Escalante-Gonz\'alez\thanks{rodolfo.escalante@ipicyt.edu.mx} }
\author{E.Campos-Cant\'on\thanks{eric.campos@ipicyt.edu.mx (Corresponding author)}}
\affil{Divisi\'on de Matem\'aticas Aplicadas, Instituto Potosino de Investigaci\'on Cient\'ifica y Tecnol\'ogica A. C., Camino a la Presa San Jos\'e 2055, Col. Lomas 4 Secci\'on, C.P. 78216, San Luis Potos\'i, S.L.P., M\'exico.}
\date{}

\begin{document}
	\maketitle
\begin{abstract}
	In this paper an approach to generate hidden attractors based on piecewise linear (PWL) systems is studied. The approach consists of the coexistence of self-excited attrators and hidden attractors, {\it i.e.}, the equilibria of the system are immersed in the basin of attraction of the seft-excited attractors, so hidden attractors appear around these  self-excited attractors. 
	
	The approach starts by generating a double-scroll attractor based on two equilibria presented heteroclinic orbits.  Then, two equilibria are added to the system, so biestability is generated by displaying two self-excited attractor. In this paper we show that hidden attractors arise as a consequence of the rupture of trajectories that resemble heteroclinic orbits at a larger scale. Therefore, multistability appears naturally as an interesting phenomenon present in this class of dynamical systems via  hidden attractors and self-excited attractors. 
	
	The study suggests the feasibility of the geometric design of new classes of multistable systems with coexistence of the two classes of attractors or even multistable systems without equilibria at all.
	
\end{abstract}

%%%%%%%%%%%%%
%%Section 1%%
%%%%%%%%%%%%%
\section{Introduction}
In the area of nonlinear dynamical systems, the  analysis of multistability is an interesting feature to study. Multistability is usually related to the existence of more than one attractor, different scenarios of multistability are reported in \cite{Anzo18}.

On the other hand, there are two classes of attractors according to \cite{Leonov11}, which are defined as follows: the first class is given by those classical attractors excited from unstable equilibria called {\it self-excited attractors} whose basin of attraction intersects with an arbitrarily small open neighborhood of equilibria, \cite{Dudkowski16} and they are not difficult to find via numerical methods, and the second class is called {\it hidden attractors} whose basin of attraction does not contain neighborhoods of equilibria. The localization of this last class represents a more difficult task which has led to interesting approaches as the analytical-numerical algorithm suggested in \cite{Leonov11} for the localization of hidden attractors of Chua's circuit.

Hidden attractors can be found in coexistence with self-excited attractors or even in systems without equilibria whose study has been the focus of recent works.

Arnold Sommerfeld worked with one of the first dynamical systems with oscillating behavior but no equilibria \cite{Kiseleva16}. In 1994 a conservative system without equilibria which presents a chaotic flow was reported in \cite{Sprott94}. This system known as Sprott case A presents two quadratic nonlinearities and it is a particular case of the Nose-Hoover system \cite{Hoover95}. After this work, several three dimensional systems without equilibria with chaotic attractors have been reported, as the one in \cite{Wei11} with two quadratic nonlinearities based on the systems Sprott case D, the one in \cite{Wang13} with three quadratic nonlinearities or the piecewise linear system reported in \cite{Escalante-Gonzalez17-1}.
In \cite{Sajad13} three methods are used to produce seventeen three dimensional systems without equilibria with chaotic flows which present only quadratic nonlinearities. 

Four dimensional systems without equilibria with chaotic or hyperchaotic attractors have also been reported, for instance, systems with quadratic and cubic nonlinearities with hyperchaotic attractors are reported in \cite{Wang12-2} and \cite{Pham16}.
The first piecewise linear system without equilibria which exhibits a hyperchaotic attractor is reported in \cite{Li14}, it is the result of the approximation made to the quadratic nonlinearities of an extended diffusionless Lorenz system. 

In \cite{Tahir15} a four dimensional system without equilibria with chaotic multiwing butterfly attractors is presented.

Since the double scroll attractor in Chua's circuit there exists an interest to generate double scroll and multiscroll attractors. 

Some approaches for self-excited attractors have been reported in \cite{Suykens97,Tang01,Yalcin01,Campos-Canton15,Campos-Canton16}.
Recently, in \cite{Escalante-Gonzalez17-2} an approach for the generation of multiscroll hidden attractors with any number of scrolls in a system without equilibria was introduced. In \cite{Hu16} two systems with multiscroll hidden attractors are constructed by introducing nonlinear functions into Sprott A system. In \cite{Jafari16} a no-equilibrium system with multiscroll hidden chaotic sea is introduced. In \cite{Hu17} a memristive chaotic system is proposed, in which multi-scroll hidden attractors and multi-wing hidden attractors which are sensitive to the transient simulation can be observed.

In \cite{Ontanon17} a study on the widening of the basins of attraction in piecewise linear multistable systems is presented, also, a system with two double scroll self excited attractors and one double scroll hidden attractor is introduced.
Based on this result it is natural to think in the possibility of generate hidden attractors via multistable systems with double scroll self-excited attractors.

In this work an approach to generate hidden attractors via piecewise linear systems (PWL) with double scroll self-excited attractors is studied. The study reveals a relation between the emergence of a hidden attractor and the existence of trajectories  that resemble heteroclinic orbits at a larger scale joining the self-excited attractors. This suggests the feasibility of the geometric design of new classes of multistable systems.

The structure of the article is as follows: In Section \ref{sec:heteroclinic} a class of piecewise linear systems with double scroll self-excited chaotic attractors is introduced; In Section \ref{Sec:multicroll_attractors} additional equilibria is considered to generate two self-excited attractors; In Section \ref{Sec:SEA} the transitory behavior of the trajectories around the self-excited attractors of the system is studied; In Section \ref{Sec_HA} the relation between the emergence of a hidden attractor and the existence of trajectories  that resemble heteroclinic orbits at a larger scale joining the self-excited attractors is discussed; Finally, conclusions are given in section \ref{Sec:conclusion}.

%%%%%%%%%%%%%
%%Section 2%%
%%%%%%%%%%%%%
\section{Heteroclinic chaos}\label{sec:heteroclinic}
To introduce the approach, let us first consider a partition $P$ of the metric space $X\subset\R^3$, endowed with the Euclidean metric $d$. Let $P=\{P_1,\ldots,P_\eta\}$ $(\eta>1)$ be a finite partition of $X$, that is, $X=\bigcup_{1\leq i\leq \eta}P_i$, and $P_i\cap P_j=\emptyset$ for $i\neq j$. Each element of the set $P$ is called an atom and each atom contains a saddle equilibrium point. Due to these atoms $P_i$ have a saddle equilibrium point then there is a stable manifold and another unstable manifold inside each atom. These stable manifolds $W^s$ and unstable manifolds $W^u$  are necessary for the mechanism of expansion and contraction present in chaotic dynamics. 

Let $T:X\to X$, with $X\subset \R^3$, be a piecewise linear dynamical system whose dynamics is given by a family of sub-systems of the form 
\begin{equation}\label{eq:affine}
	\dot{\bf x}=A{\bf x}+f({\bf x})B,
\end{equation}
where ${\bf x}=(x_1,x_2,x_3)^T\in\mathbb{R}^3$ is the state vector, and  $A=\{\alpha_{ij}\} \in \mathbb{R}^{3 \times 3}$ is a linear operator, $B=(\beta_{ 1}, \beta_{ 2}, \beta_{3})^T$  is a constant vector, and $f$ is a functional. The vector $f({\bf x})B$ is a constant vector in each atom $P_i$ such that the equilibria is given by ${\bf x}^*_{eq_i}=(x^*_{{1}_{eq_i}},x^*_{{2}_{eq_i}},x^*_{{3}_{eq_i}})^T=-f({\bf x})A^{-1}B\in P_i$, with $i=1,\ldots,\eta$.

To induce oscillations of the flow around the equilibria ${\bf x}^*_{eq_i}$ it is necessary to assign a negative real eigenvalue $\lambda_1=c$ to matrix $A$ with the corresponding eigenvector $v_{1}$, and a pair of complex conjugate eigenvalues with positive real part $\lambda_2=a+ib$ and $\lambda_3=a-ib$ with the corresponding eigenvectors $v_{2}$ and $v_{3}$. Thus the stable and unstable manifolds are given by $W^s_{{\bf x}^*_{eq_i}}=\{{\bf x}+{\bf x}^*_{eq_i}:{\bf x}\in span\{v_{1}\}  \}$ and  $W^u_{{\bf x}^*_{eq_i}}=\{{\bf x}+{\bf x}^*_{eq_i}:{\bf x}\in span\{v_{2},v_{3}\}  \}$. These sets generate the contraction and expansion of the trajectories. 

The matrix of the linear operator $A$ is defined as follows:

\begin{equation}
	A=\begin{pmatrix}\label{eq:Amatrix}
		\frac{a}{3} + \frac{2c}{3}&  b& \frac{2c}{3} - \frac{2a}{3}\\
		-\frac{b}{3}&  a&           \frac{2b}{3}\\
		\frac{c}{3} - \frac{a}{3}& -b&     \frac{2a}{3} + \frac{c}{3}\\
	\end{pmatrix},
\end{equation}
whose eigenvectors are
\begin{equation}
	v_{1}=\begin{pmatrix}
		1\\0\\\frac{1}{2}
	\end{pmatrix}, \
	v_{2}=\begin{pmatrix}
		0\\-1\\0
	\end{pmatrix}, \ 
	v_{3}=\begin{pmatrix}
		-1\\0\\1
	\end{pmatrix}.
\end{equation}
\color{black}
In this work we denote the local stable and unstable manifolds of an equilibrium point ${\bf x}^*_{eq}$ as $W^s_{{\bf x}^*_{eq}}$ and $W^u_{{\bf x}^*_{eq}}$, respectively, and they are responsible for connecting the equilibria of a dynamical system. Recall that a path in phase space which joins two different equilibrium points is called a {\it heteroclinic orbit}. And a path that starts and ends at the same point is called a {\it homoclinic orbit}.

We also denoted the closure of a set $P_i$ as $cl(P_i)$. Thus, for each pair of atoms $P_i$ and $P_j$, $i\neq j$, if $cl(P_i)\cap cl(P_j)\neq \emptyset$ then these atoms are adjacent and the switching surface between them is given by the intersection, {\it i.e.},  $SW_{ij}=cl(P_i)\cap cl(P_j)$.

\begin{ass}
	The divergence of the PWL system \eqref{eq:affine} considering the linear operator $A$ given by \eqref{eq:Amatrix} is $\nabla=2a+c$, so the system is dissipative in each atom of the partition $P$ if $2a<|c|$.  
\end{ass}

With the atoms of a  $P$ partition containing a saddle equilibrium point in each of them as defined above, it is possible to generate heteroclinic orbits. At least two equilibria are necessary to  generate a heteroclinic orbit, so we start by considering a partition with two atoms $P=\{P_{1},P_{2}\}$, the constant vector $B\in\mathbb{R}^3$ is defines as follows:
\begin{equation}\label{eq1:vectorB}
	B=\begin{pmatrix}
		-\frac{a}{3} - \frac{2c}{3}\\
		\frac{b}{3}\\
		\frac{a}{3} -\frac{c}{3}\\
	\end{pmatrix},
\end{equation}
and the functional $f$ is given by
\begin{equation}\label{eq1:functional}
	f({\bf x})=\left\{ 
	\begin{array}{rl}
		-\alpha, & {\bf x}\in P_1;\\
		\alpha, & {\bf x}\in P_2;
	\end{array}\right .
\end{equation}
with $\alpha>0$. So the equilibria are at ${\bf x}^*_{eq_1}=(-\alpha,0,0)^T\in P_{1}$ and ${\bf x}^*_{eq_2}=(\alpha,0,0)^T\in P_{2}$, and the stable and the unstable manifolds are given by 
\begin{align*}
	& W^s_{{\bf x}^*_{eq_1}}=\left\{ {\bf x}\in \mathbb{R}^3| x_1+\alpha=2x_3,x_2=0\right\},\\
	& W^u_{{\bf x}^*_{eq_1}}=\left\{ {\bf x}\in \mathbb{R}^3|  x_1+x_3=-\alpha \right\},\\
	& W^s_{{\bf x}^*_{eq_2}}=\left\{ {\bf x}\in \mathbb{R}^3|  x_1-\alpha=2x_3,x_2=0\right\},\\
	& W^u_{{\bf x}^*_{eq_2}}=\left\{ {\bf x}\in \mathbb{R}^3|  x_1+x_3=\alpha \right\}.\\
\end{align*}

The next proposition gives the necessary and sufficient conditions to generate heteroclinic connections between equilibria of a  piecewise linear dynamical system. Thus, the appearance of a chaotic attractor in this type of systems is possible.

\begin{proposition}\label{pro:Heteroclinic_orbit}
	The hyperbolic system given by \eqref{eq:affine}, \eqref{eq:Amatrix}, \eqref{eq1:vectorB} and \eqref{eq1:functional} generates a pair of heteroclinic orbits if the switching surface between the atoms $P_1$ and $P_2$ is given by the plane $SW=cl(P_1)\cap cl(P_2)$ with $\{{\bf x}\in\mathbb{R}^3:x_3\geq0\}\cap SW \in P_1$ and $\{{\bf x}\in\mathbb{R}^3:x_3< 0\}\cap SW \in P_2$ and also pass through the intersection points $cl(W^s_{{\bf x}^*_{eq_1}}) \cap cl(W^u_{{\bf x}^*_{eq_2}})\neq \emptyset$, and $cl(W^s_{{\bf x}^*_{eq_2}}) \cap cl(W^u_{{\bf x}^*_{eq_1}})\neq \emptyset$ and the line: $x_1=(x^*_{{1}_{eq_1}}+x^*_{{1}_{eq_2}})/2, x_3=0$. 
\end{proposition}

{\it Proof:} We want to show that there exist initial conditions ${\bf x}_{01}, {\bf x}_{02}\in SW$, such that two solution curves $\varphi({\bf x}_{01},t)$ and $\varphi({\bf x}_{02},t)$ of the hyperbolic system given by \eqref{eq:affine}, \eqref{eq:Amatrix}, \eqref{eq1:vectorB} and \eqref{eq1:functional} fulfill that $\varphi({\bf x}_{01},t)\to {\bf x}^*_{eq_1}$ and $\varphi({\bf x}_{02},t)\to {\bf x}^*_{eq_2}$ as $t\to \infty$ and $\varphi({\bf x}_{01},t)\to {\bf x}^*_{eq_2}$ and $\varphi({\bf x}_{02},t)\to {\bf x}^*_{eq_1}$ as $t\to -\infty$, in particular, these initial conditions correspond to the intersection points $cl(W^s_{{\bf x}^*_{eq_1}}) \cap cl(W^u_{{\bf x}^*_{eq_2}})$, and $cl(W^s_{{\bf x}^*_{eq_2}}) \cap cl(W^u_{{\bf x}^*_{eq_1}})$. 

From \eqref{eq:Amatrix} the linear operator $A$ can be expressed as
\begin{equation}
	A=QEQ^{-1},
\end{equation}
where $Q=[v_1 \; v_2 \; v_3]$ and 
\begin{equation}
	E=\begin{pmatrix}
		c&0&0\\
		0&a&-b\\
		0&b&a\\
	\end{pmatrix}.
\end{equation}

According to the stable and unstable manifolds, the intersection points  are give by

$${\bf x}_{in_1}=cl(W^s_{{\bf x}^*_{eq_1}}) \cap cl(W^u_{{\bf x}^*_{eq_2}})=\left(\frac{\alpha}{3},0,\frac{2\alpha}{3}\right)^T.$$	
$${\bf x}_{in_2}=cl(W^s_{{\bf x}^*_{eq_2}}) \cap cl(W^u_{{\bf x}^*_{eq_1}})=\left(-\frac{\alpha}{3},0,-\frac{2\alpha}{3}\right)^T,$$

These points ${\bf x}_{in_1}$ and ${\bf x}_{in_2}$ belong to $SW$ and  ${\bf x}_{in_1}\in P_1$ and ${\bf x}_{in_2}=P_2$. Because these points ${\bf x}_{in_1}$ and ${\bf x}_{in_2}$  belong to the stable manifolds  $W^s_{{\bf x}^*_{eq_1}}$ and $W^s_{{\bf x}^*_{eq_2}}$, respectively, they are points whose trajectories remain in atoms $P_1$ and $P_2$, respectively.

By definition ${\bf x}^*_{1_{eq_1}}=-{\bf x}^*_{1_{eq_2}}$, then the $x_2$ axis belongs to the plane $SW$. The sets $cl(W^u_{{\bf x}^*_{eq_i}})\cap SW$, for $i=1,2$, can be written as: 

\begin{equation}\label{eq:interseccion}
	\{(0,\epsilon,0)^T+{\bf x}_{in_i}:\epsilon\in\mathbb{R}\}, \text{ for } i=1,2.
\end{equation}

Consider the following changes of coordinates ${\bf z}^{(i)}=Q^{-1}({\bf x}-{\bf x}^*_{eq_i})$, for $i=1,2$. Then the vector field in ${\bf z}^{(i)}$ coordinates for the space given by the atom $P_i$ is given by $\dot{{\bf z}}^{(i)}=E{\bf z}^{(i)}$, with $i=1,2$.

Since $Q^{-1}(0,\epsilon,0)^T=(0,-\epsilon,0)^T$, the sets given by \eqref{eq:interseccion} in ${\bf z}^{(i)}$ coordinates are given as follows
\begin{equation}\label{eq:interseccionz}
	\{(0,\epsilon,0)^T+Q^{-1}({\bf x}_{in_i}-{\bf x}^*_{eq_i}):\epsilon\in\mathbb{R}\}, \text{ for } i=1,2,
\end{equation}

where ${\bf z}^{(i)}_{in_i}=Q^{-1}({\bf x}_{in_i}-{\bf x}^*_{eq_i})= \left((-1)^{i+1}4\alpha/3,0,0\right)^T$ is a point on the ${ z}_1^{(i)}$ axis that corresponds to the transformation of the intersection points ${\bf x}_{in_i}\in P_i$ to ${\bf z}^{(i)}_{in_i}\in \{ Q^{-1}({\bf x}-{\bf x}^*_{eq_i}): {\bf x}\in P_i \}$, for $i=1,2$.

When $t>0$, $\varphi({\bf x}_{in_1},t)$ remains in the atom $P_1$, the transformation of ${\bf x}_{in_1}$  under $Q^{-1}({\bf x}_{in_1}-{\bf x}^*_{eq_1})$ is ${\bf z}^{(1)}_{in_1}=( 4\alpha/3,0, 0)^T$. In a similar way, when $t>0$, $\varphi({\bf x}_{in_2},t)$,  remains in the atom $P_2$, the transformation of ${\bf x}_{in_2}$ under $Q^{-1}({\bf x}_{in_2}-{\bf x}^*_{eq_2})$ is ${\bf z}^{(2)}_{in_2}=(-4\alpha/3,0, 0)^T$. So  ${\bf z}^{(i)}_{in_i}$ belongs to the stable manifold $W^s_{{\bf z}^*_{eq_i}}$, for $i=1,2$, then the trajectories ${\bf z}^{(i)}(t)=e^{Et}{\bf z}^{(i)}_{in_i} \to 0$ when $t\to \infty$. This implies that
\[\lim\limits_{t\to\infty}\varphi({\bf x}_{in_1},t)={\bf x}^*_{eq_1}, \ \text{and} \ \lim\limits_{t\to \infty}\varphi({\bf x}_{in_2},t)={\bf x}^*_{eq_2}. \]

When $t<0$, $\varphi({\bf x}_{in_1},t)$ leaves the atom $P_1$ and enters to atom $P_2$, the transformation of ${\bf x}_{in_1}$ under $Q^{-1}({\bf x}_{in_1}-{\bf x}^*_{eq_2})$ to ${\bf z}^{(2)}_{in_1}\in \{ Q^{-1}({\bf x}-{\bf x}^*_{eq_2}): {\bf x}\in cl(P_2) \}$ is ${\bf z}^{(2)}_{in_1}=(0,0, 2\alpha/3)^T$. In a similar way, when $t<0$, $\varphi({\bf x}_{in_2},t)$ leaves the atom $P_2$ and enters to atom $P_1$, the transformation of ${\bf x}_{in_2}$ under $Q^{-1}({\bf x}_{in_2}-{\bf x}^*_{eq_1})$ to ${\bf z}^{(1)}_{in_2}\in \{ Q^{-1}({\bf x}-{\bf x}^*_{eq_1}): {\bf x}\in cl(P_1) \}$ is ${\bf z}^{(1)}_{in_2}=(0,0,-2\alpha/3)^T$. Thus, ${\bf z}^{(j)}_{in_i}=(0,0,(-1)^{j}2\alpha/3)^T$ is a point on the axis $z_3^{(j)}$ and belongs to $cl(W^u_{{\bf z}^*_{eq_j}})$ for $i,j\in\{1,2\}$ and $i\neq j$ .

With the uncoupled system  in ${\bf z}^{(i)}$ coordinates, we can analyze the flow on the plane $z_2^{(i)}-z_3^{(i)}$ and see how the flow converges at the equilibrium point ${\bf z}^{*(i)}_{eq_j}$ when $t\to -\infty$.
\begin{equation}
	\begin{array}{l}
		\dot{z_2}^{(i)}=az_2^{(i)}-bz_3^{(i)},\\
		\dot{z_3}^{(i)}=bz_2^{(i)}+az_3^{(i)},
	\end{array}
\end{equation}
\begin{equation}
	z_2^{(i)}\dot{z_2}^{(i)}+z_3^{(i)}\dot{z_3}^{(i)}=a\left(\left(z_2^{(i)}\right)^2+\left(z_3^{(i)}\right)^2\right),
\end{equation}
if $r^2=\left(z_2^{(i)}\right)^2+\left(z_3^{(i)}\right)^2$ then $r\dot{r}=ar^2$
\begin{equation}
	\dot{r}=ar,
\end{equation}
\begin{equation}\label{eq:polarsolution}
	r=r_0e^{at}.
\end{equation}

As $0<a\in \mathbb{R}$, so $r\to 0$ when $t\to -\infty$.  Then the trajectories ${\bf z}^{(i)}(t)=e^{Et}{\bf z}^{(i)}_{in_i} \to 0$ when $t\to -\infty$.
This implies that 
\[\lim\limits_{t\to-\infty}\varphi({\bf x}_{in_1},t)={\bf x}^*_{eq_2}, \ \text{and} \ \lim\limits_{t\to-\infty}\varphi({\bf x}_{in_2},t)={\bf x}^*_{eq_1} ,\]

Thus the heteroclinic orbits are defined as
\[HO_1= \{{\bf x}\in \varphi({\bf x}_{in_1},t):t\in(-\infty,\infty)\},\]
\[HO_2= \{{\bf x}\in \varphi({\bf x}_{in_2},t):t\in(-\infty,\infty)\}.\]

\begin{flushright}
	$\square$
\end{flushright}

For the system given by \eqref{eq:affine}, \eqref{eq:Amatrix}, \eqref{eq1:vectorB} and \eqref{eq1:functional} it is possible to find several points ${\bf x}_0\in HO_i$ such that $|{\bf x}_{eq_i}-{\bf x}_0|<\epsilon$ with $\epsilon$ arbitrarily small and $i=1,2$.

For example, if an arbitrary initial condition ${\bf z}^{(i)}_0=(z^{(i)}_{10},z^{(i)}_{20},z^{(i)}_{30})^T$ belongs to $z^{(i)}_3$ axis and the heteroclinic orbit $HO_i$, then the intersection point $z_{in_i}$ is reached at time $t_f=\frac{2k\pi}{b}$, with $k\in\mathbb{Z}^+$. To see this, we consider the transformation introduced before ${\bf z}^{(i)}=Q^{-1}({\bf x}-{\bf x}^*_{eq_i})$.
Then from the solution in ${\bf z}^{(i)}$ coordinates in $cl(P_i)$: 
\begin{equation}
	\begin{array}{l}
		z_1(t)^{(i)}=z_{10}^{(i)}e^{c t},\\
		z_2(t)^{(i)}=z_{20}^{(i)}e^{a t}\cos(b t) -z_{30}^{(i)}e^{a  t}\sin(b t),\\
		z_3(t)^{(i)}=z_{20}^{(i)}e^{a t}\sin(b t) +z_{30}^{(i)}e^{a  t}\cos(b t).
	\end{array}
\end{equation}
Let us assume an initial condition of the form ${\bf z}_0^{(i)}=(0,0,z^{(i)}_{30})^T$ due to it is in the $z^{(i)}_3$ axis, then 
\begin{equation}\label{Sol_ini_cond_ejeZ}
	\begin{array}{l}
		z_1(t)^{(i)}=0,\\
		z_2(t)^{(i)}= -z_{30}^{(i)}e^{a  t}\sin(b t),\\
		z_3(t)^{(i)}=z_{30}^{(i)}e^{a  t}\cos(b t).
	\end{array}
\end{equation}
Because the initial condition belongs to the heteroclinic orbit $HO_i$, then from the second equation of \eqref{Sol_ini_cond_ejeZ} the trajectory  $\varphi ({\bf z}^{(i)}_0,t)$ reaches the switching surface at the intersection point ${\bf z}^{(i)}_{in_j}=(z^{(i)}_{1in_j},z^{(i)}_{2in_j},z^{(i)}_{3in_j})=\left(0,0,(-1)^{i}2\alpha/3\right)^T$ at time $t_f=\frac{2k\pi}{b}$, with $k\in\mathbb{Z}^+$.
Thus
\begin{equation}
	z^{(i)}_{3in_j}=z_{30}^{(i)}e^{  \frac{2ka\pi}{b}},
\end{equation}
\begin{equation}
	z_{30}^{(i)}=z^{(i)}_{3in_j}e^{-\frac{2ka\pi}{b}}.
\end{equation}

Then the initial condition ${\bf x}^i_0=Q{\bf z}^{(i)}_0+{\bf x}^*_{eq_i}$, for $P_1$
\begin{equation}
	{\bf x}^1_{0}=\begin{pmatrix}\label{eq:x01}
		\frac{2}{3}\alpha e^{-\frac{2ka\pi}{b}}-\alpha\\0\\-\frac{2}{3}\alpha e^{-\frac{2ka\pi}{b}}
	\end{pmatrix},
\end{equation}
and $P_2$
\begin{equation}\label{eq:x02}
	{\bf x}^2_{0}=\begin{pmatrix}
		-\frac{2}{3}\alpha e^{-\frac{2ka\pi}{b}}+\alpha\\0\\\frac{2}{3}\alpha e^{-\frac{2ka\pi}{b}}
	\end{pmatrix}.
\end{equation}
Then we could also express the heteroclinic orbits as
\[HO_1= \{{\bf x}\in \varphi({\bf x}^1_{0},t):t\in(-\infty,0]\} \cup\{{\bf x}\in \varphi({\bf x}^1_{0},t):t\in[0,(2k\pi)/b), k\in\mathbb{Z}^+\} \cup \{{\bf x}\in \varphi({\bf x}_{in_1},t):t\in[0,\infty)\},\]
\[HO_2= \{{\bf x}\in \varphi({\bf x}^2_{0},t):t\in(-\infty,0]\}\cup \{{\bf x}\in \varphi({\bf x}^2_{0},t):t\in[0,(2k\pi)/b),k\in\mathbb{Z}^+\}\cup \{{\bf x}\in \varphi({\bf x}_{in_2},t):t\in[0,\infty)\}.\]

The trajectory ${\bf x}(t)$ of the PWL system  can be calculated by ${\bf x}^i (t)=e^{At} {\bf x}^i_0$ in each atom $P_i$, where ${\bf x}^i= {\bf x} +{\bf x}^*_{eq_i}$ and ${\bf x}^i_0$ is the initial condition when the trajectory enter to the atom $P_i$,  $i=1,2$. Then
$${\bf x}^i(t)=QE(t)Q^{-1} {\bf x}^i(0),$$
where $Q$ is the invertible matrix defined by the eigenvectors of $A$ and
\begin{equation}
	E(t)= \begin{pmatrix}
		e^{c t} & 0 & 0 \\
		0 & e^{a t}\cos(b t) & -e^{a t}\sin(b t) \\
		0 & e^{a t}\sin(b t) & e^{a  t}\cos(b t) \end{pmatrix} .
\end{equation} 
So the exact solution is given as follows: 

\begin{equation}
	\begin{array}{llll}
		x_1^i(t)=&x^i_{10}(e^{at}\cos(bt) + 2e^{ct})/3&-x^i_{20}e^{at}\sin(bt) &+x^i_{30}(2e^{ct}-2e^{at}\cos(bt))/3, \\
		x_2^i(t)=&x^i_{10}(e^{at}\sin(bt))/3 &+x^i_{20}e^{at}\cos(bt) & -x^i_{30}(2e^{at}\sin(bt))/3,\\
		x_3^i(t)=&x^i_{10}(e^{ct}-e^{at}\cos(bt))/3&+e^{at}x^i_{20}\sin(bt)&+x^i_{30}(2e^{at}\cos(bt)+e^{ct})/3. 
	\end{array}
\end{equation}

The flow of the system $\varphi({\bf x}_0)$ is dissipative in each atom for all initial condition ${\bf x}_0\in P_i - W^u_{{\bf x}_{eq_i}} \subset X$. If the initial condition ${\bf x}_0\in W^u_{{\bf x}_{eq_i}}$ then ${\bf x}_0$ is a linear combination of $v_{2}$ and $v_{3}$, {\it e.g.} ${\bf x}_0=(x_{10},x_{20},-x_{10})^T$.

\begin{equation}
	\begin{array}{lr}
		x_1^i(t)=& e^{at}(x^i_{10}\cos(bt) - x^i_{20}\sin(bt)),\\
		x_2^i(t)=&  e^{at}(x^i_{20}\cos(bt) + x^i_{10}\sin(bt)),\\
		x_3^i(t)=& -e^{at}(x^i_{10}\cos(bt) - x^i_{20}\sin(bt)).
	\end{array}
\end{equation}

\begin{ass}\label{as:ab}
	The oscillations around the equilibrium point $x^*_{eq_i}$ depend on parameters $a$ and $b$, we consider $b/a>10$.
\end{ass}

\begin{figure}[!ht]
	\centering	
	\subfloat[]{\includegraphics[width=0.5\columnwidth]{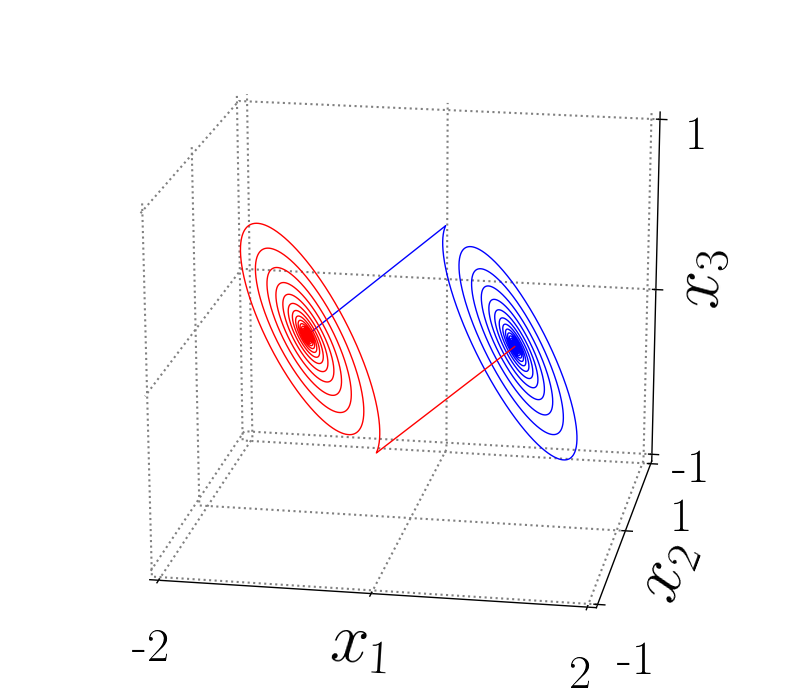}\label{fig:heteroclinic}} \hfill
	\subfloat[]{\includegraphics[width=0.5\columnwidth]{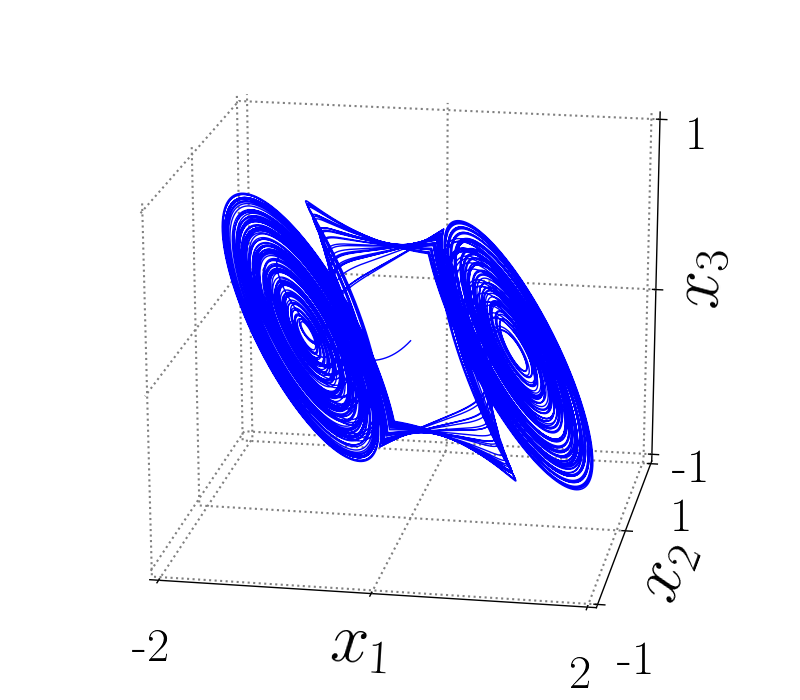}\label{fig:doble}}
	\caption{\label{fig:dobleorbits} In (a) the heteroclinic loop of the system \eqref{eq:affine}, \eqref{eq:Amatrix},\eqref{eq1:vectorB} and \eqref{eq1:functional} with the switching surface $\{ x\in\mathbb{R}^3:2x_1-x_3=0\}$ and the parameters $a=0.2,b=5,c=-3,\alpha=1$, and in (b) a double-scroll attractor that emerges from an heteroclinic orbit using the following initial condition ${\bf x}_0=(0,0,0)^T$.}
\end{figure}

As example, we consider the system \eqref{eq:affine}, \eqref{eq:Amatrix},\eqref{eq1:vectorB} and \eqref{eq1:functional} with the switching surface and parameters  $a,\ b,\ c,\ \alpha$ as follows.
\begin{example}
	$SW=\{ x\in\mathbb{R}^3:2x_1-x_3=0\}$,  with $\{x\in\mathbb{R}^3:x_3\geq 0\}\cap SW \in P_1$, and $\{x\in\mathbb{R}^3:x_3< 0\}\cap SW \in P_2$ and the parameters $a=0.2,b=5,c=-3,\alpha=1$.
\end{example}
The  above defined system fulfills the proposition \ref{pro:Heteroclinic_orbit}, so it presents a heteroclinic orbit. From \eqref{eq:x01} and \eqref{eq:x02} two initial conditions \\
${\bf x}_{01}=(-0.9999976751050959, 0, -2.3248949041393315e-6)^T$ and\\
${\bf x}_{02}=(0.9999976751050959, 0, 2.3248949041393315e-6)^T$ are chosen with $k=50$ to simulate the two heteroclinic orbits shown in the Figure~\ref{fig:heteroclinic}. Then heteroclinic chaos emerges from this system, in particular, a double scroll attractor is generated as it is shown in the Figure~\ref{fig:doble}, for the following initial condition ${\bf x}_0=(0,0,0)^T$.

The unstable manifolds $W^u_{{\bf x}^*_{eq_1}}=\{{\bf x}\in \R^3 : x_1+x_3+1=0\}$ and $W^u_{{\bf x}^*{eq_2}}=\{{\bf x}\in \R^3 : x_1+x_3-1=0\}$ and the stable manifolds $W^s_{{\bf x}^*_{eq_1}}=\{{\bf x}\in \R^3 : \frac{x_1+1}{2}=x_3; x_2=0\}$ and $W^s_{{\bf x}^*_{eq_2}}=\{{\bf x}\in \R^3 : \frac{x_1-1}{2}=x_3; x_2=0\}$.
The intersection points are given by $cl(W^s_{{\bf x}^*_{eq_2}})\cap cl(W^u_{{\bf x}^*_{eq_1}})=(-\frac{1}{3},0,-\frac{2}{3})^T$,  $cl(W^s_{{\bf x}^*_{eq_1}}) \cap cl(W^u_{{\bf x}^*_{eq_2}})=$ $(\frac{1}{3},0,\frac{2}{3})^T$.

\begin{proposition}\label{prop:several_HO}
	If the partition $P$ contains more that two atoms $\{P_1,P_2,\ldots,P_k\}$, with $2<k\in \mathbb{Z}^+$, and each atom is a hyperbolic set defined as above. Furthermore, the atoms by pairs $P_i$ and $P_{i+1}$ fulfill the proposition  \ref{pro:Heteroclinic_orbit}. Then the system generates $2(k-1)$ heteroclinic orbits.
\end{proposition}
{\it Proof:} A direct consequence of the proposition  \ref{pro:Heteroclinic_orbit} \begin{flushright}
	$\square$
\end{flushright}

%%%%%%%%%%%%%
%%Section 3%%
%%%%%%%%%%%%%
\section{Emergence of multiscroll attractors through multiple heteroclinic orbits}\label{Sec:multicroll_attractors}
According to the proposition \ref{prop:several_HO}, it is possible to generate multiscrol attractors based on multiple heteroclinic orbits. So in this Section \ref{Sec:multicroll_attractors} we consider more than two hyperbolic sets of a partition, without loss of generality,  the following partition is considered $P=\{P_1,P_2, P_3,P_4\}$, to generate a multiscroll attractor.

Consider the  piecewise linear dynamical system \eqref{eq:affine}, with $A$ and $B$ given by \eqref{eq:Amatrix} and \eqref{eq1:vectorB}, respectively, thus  the functional $f({\bf x})$ is defined in the four atoms as follows:

\begin{equation}\label{eq:funcion4atomos}
	f({\bf x})=\left\{\begin{array}{ll}
		-\alpha-\gamma,&{\bf x}\in P_1;\\
		\alpha-\gamma,&{\bf x}\in P_2;\\
		-\alpha+\gamma,&{\bf x}\in P_3;\\
		\alpha+\gamma,&{\bf x}\in P_4;
	\end{array}\right.
\end{equation}
where $\alpha,\gamma>0$. The equilibria are at:
\begin{equation}
	{\bf x}^*_{eq_1}=\begin{bmatrix}
		-(\gamma+\alpha)\\0\\0
	\end{bmatrix},\ 
	{\bf x}^*_{eq_2}=\begin{bmatrix}
		-(\gamma-\alpha)\\0\\0
	\end{bmatrix},\
	{\bf x}^*_{eq_3}=\begin{bmatrix}
		(\gamma-\alpha)\\0\\0
	\end{bmatrix},\
	{\bf x}^*_{eq_4}=\begin{bmatrix}
		(\gamma+\alpha)\\0\\0
	\end{bmatrix},
\end{equation}
so ${\bf x}^*_{eq_1}\in P_1$, ${\bf x}^*_{eq_2}\in P_2$, ${\bf x}^*_{eq_3}\in P_3$ and ${\bf x}^*_{eq_4}\in P_4$. The location of the equilibria according to the parameters $0<\alpha$ and $0<\gamma$ is as follows:
\begin{itemize}
	\item The equilibria are on the $x_1$ axis and for $\alpha=\gamma$ the system only have three equilibria, otherwise it has four equilibria. 
	\item For $\alpha <\gamma$ the distance of the equilibria ${\bf x}^*_{eq_1}$ and ${\bf x}^*_{eq_4}$ to the origin $O=(0,0,0)^T$ are the same $d({\bf x}^*_{eq_1},O)=d({\bf x}^*_{eq_4},O)$ and also for $d({\bf x}^*_{eq_2},O)=d({\bf x}^*_{eq_3},O)$.
	\item For $\gamma=2\alpha$, all equilibria are at the same distance $d({\bf x}^{*}_{eq_1},{\bf x}^{*}_{eq_2})=d({\bf x}^{*}_{eq_2},{\bf x}^{*}_{eq_3})=d({\bf x}^{*}_{eq_3},{\bf x}^{*}_{eq_4})=2\alpha$. 
	\item  The another case is when  $\gamma\neq 2\alpha$, and $d({\bf x}^{*}_{eq_1}, {\bf x}^{*}_{eq_2})=d({\bf x}^{*}_{eq_3},{\bf x}^{*}_{eq_4})=2\alpha$, but $d({\bf x}^{*}_{eq_2},{\bf x}^{*}_{eq_3})\neq 2\alpha$.
\end{itemize}
In this Section \ref{Sec:multicroll_attractors}, we are also interested in the case of $\gamma\neq 2\alpha$ such that $\gamma>\alpha$.

The switching surfaces are given by:
\begin{equation}\label{eq:surfaces4atoms}
	\begin{array}{l}
		SW_{12}=cl(P_1)\cap cl(P_2)=\{{\bf x}\in\mathbb{R}^3:2x_1-x_3=-2\gamma\},\\
		SW_{23}=cl(P_2)\cap cl(P_3)=\{{\bf x}\in\mathbb{R}^3:2x_1-x_3=0\},\\
		SW_{34}=cl(P_3)\cap cl(P_4)=\{{\bf x}\in\mathbb{R}^3:2x_1-x_3=2\gamma\},
	\end{array}
\end{equation}

which fulfill that

\begin{equation}\label{eq:surfaces4atoms2}
	\begin{array}{l}
		SW_{i(i+1)}\cap \{{\bf x}\in\mathbb{R}^3:x_3>0 \}\in P_i,\\
		SW_{i(i+1)}\cap \{{\bf x}\in\mathbb{R}^3:x_3 \leq 0 \}\in P_{i+1}.
	\end{array}
\end{equation}
This way of defining the switching surfaces provokes that the intersections between them and the stable manifolds contain a point, and the intersections between them and the unstable manifolds are the empty set, {\it i.e.}, $W^u_{{\bf x}^*_{eq_1}}\cap SW_{12}=\emptyset$ and $W^s_{{\bf x}^*_{eq_1}}\cap SW_{12}\neq \emptyset$.

\begin{figure}%[h!]
	\centering
	\includegraphics[width=0.5\columnwidth]{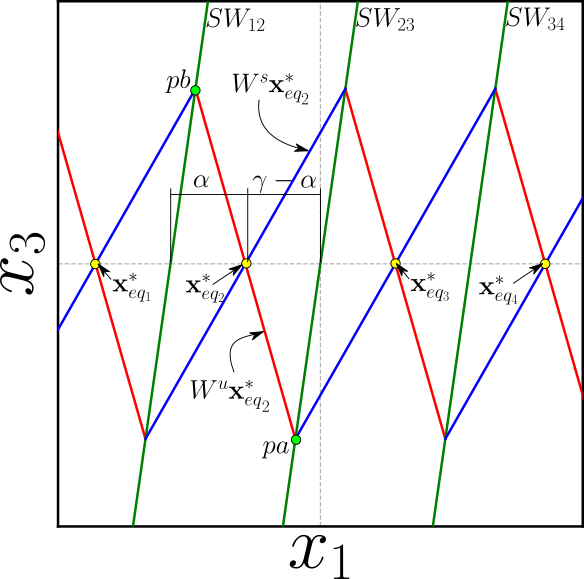}
	\caption{\label{fig:diagrama} Projection of the stable and unstable manifolds and switching planes onto the $x_1-x_3$ plane. The diagram shows the location of the unstable  manifold market with red lines, the stable  manifold market with blue lines and switching planes marked with green lines.}
\end{figure}

Let us define two points,  $pb=W^s_{{\bf x}^*_{eq_1}}\cap SW_{12}$ and $pa=W^s_{{\bf x}^*_{eq_3}}\cap SW_{23}$ as shown in the Figure~\ref{fig:diagrama}. Then $pa$ and $pb$ are given as follows

\begin{equation}
	pa=\begin{pmatrix}
		-\frac{(\gamma-\alpha)}{3}\\0\\-\frac{2(\gamma-\alpha)}{3}
	\end{pmatrix},\
	pb=\begin{pmatrix}
		\frac{\alpha}{3}-\gamma\\0\\\frac{2\alpha}{3}
	\end{pmatrix}.
\end{equation}

The set $cl(W^u_{{\bf x}^*_{eq_2}})\cap SW_{12}$ can be written as: 

\begin{equation}\label{eq:interseccionp2b}
	\{{\bf x}\in \mathbb{R}^3 :\ {\bf x}=(0,\epsilon,0)^T+pb,\epsilon\in\mathbb{R}\},
\end{equation}
and the set  $cl(W^u_{{\bf x}^*_{eq_2}})\cap SW_{23}$ can be written as: 

\begin{equation}\label{eq:interseccionp2a}
	\{{\bf x}\in \mathbb{R}^3 :\ {\bf x}=(0,\epsilon,0)^T+pa,\epsilon\in\mathbb{R}\}.
\end{equation}

Consider the transformation ${\bf z}^{(2)}=Q^{-1}({\bf x}-{\bf x}^*_{eq_2})$, since $Q^{-1}(0,\epsilon,0)^T=(0,-\epsilon,0)^T$, so ${\bf z}^{(2)}=(0,\epsilon,0)^T+Q^{-1}(pb-{\bf x}^*_{eq_2})$, where $Q^{-1}(pb-{\bf x}^*_{eq_2})$ is a point on the plane $z_2^{(2)}-z_3^{(2)}$. And ${\bf z}^{(2)}=(0,\epsilon,0)^T+Q^{-1}(pa-{\bf x}^*_{eq_2})$, where $Q^{-1}(pa-{\bf x}^*_{eq_2})$ is also a point on the plane $z_2^{(2)}-z_3^{(2)}$.
Then the set \eqref{eq:interseccionp2b} in ${\bf z}^{(2)}$ coordinates is given by
\begin{equation}\label{eq:interseccionp2bz}
	\{{\bf z}^{(2)}\in \mathbb{R}^3: {\bf z}^{(2)}=(0,\epsilon,2\alpha/3)^T,\epsilon\in\mathbb{R}\},
\end{equation}
and the set \eqref{eq:interseccionp2a} in ${\bf z}^{(2)}$ coordinates is given by
\begin{equation}\label{eq:interseccionp2az}
	\{{\bf z}^{(2)}\in \mathbb{R}^3: {\bf z}^{(2)}=(0,\epsilon,2(\alpha-\gamma)/3)^T,\epsilon\in\mathbb{R}\}.
\end{equation}

Thus the sets \eqref{eq:interseccionp2bz} and \eqref{eq:interseccionp2az} are orthogonal lines to the $z_3^{(2)}$ axis. The points $pa$ and $pb$ in ${\bf z}^{(2)}$ coordinates will be denoted as
\begin{equation}\label{eq:papbz}
	pa_z=\begin{pmatrix}
		0\\0\\\frac{2(\alpha-\gamma)}{3}
	\end{pmatrix},\
	pb_z=\begin{pmatrix}
		0\\0\\\frac{2\alpha}{3}
	\end{pmatrix}.
\end{equation}

With the uncoupled system  in ${\bf z}^{(2)}$ coordinates we can analyze the flow on the plane $z_2^{(2)}-z_3^{(2)}$ close to ${\bf z}^{*(2)}_{eq_2}$. 
\begin{equation}
	\begin{array}{l}
		\dot{z_2}^{(2)}=az_2^{(2)}-bz_3^{(2)},\\
		\dot{z_3}=bz_2^{(2)}+az_3^{(2)},
	\end{array}
\end{equation}
\begin{equation}
	r\dot{r}=z_2^{(2)}\dot{z_2}^{(2)}+z_3^{(2)}\dot{z_3}^{(2)}=ar^2,
\end{equation}
\begin{equation}
	\dot{r}=ar,
\end{equation}
\begin{equation}\label{eq:polarsolutionp2}
	r=r_0e^{at}.
\end{equation}

It follows from \eqref{eq:papbz} that  if $\alpha=\gamma-\alpha$ then the points $pa_z$ and $pb_z$ are at the same distance from ${\bf z}^{*(2)}_{eq_2}=(0,0,0)^T$. Thus, from \eqref{eq:polarsolutionp2}, it follows that the trajectories with initial conditions $pa_z$ and $pb_z$ remain in $P_2$ for all $t<0$.

Our case study is $\gamma-\alpha\neq\alpha$, such that $\gamma>\alpha$.
Let us consider the case $\gamma-\alpha>\alpha$, it can be seen from \eqref{eq:papbz} that $pb_z$ is closer to ${\bf z}^{(2)}_{eq_2}$ than $pa_z$, this is, $d(pb_z,{\bf z}^{(2)}_{eq_2})<d(pa_z,{\bf z}^{(2)}_{eq_2})$. Then, if $\gamma$ is sufficiently big with respect to $\alpha$, the trajectory with initial condition $pa_z$  will eventually reach the set   given by \eqref{eq:interseccionp2bz} for $t<0$, {\it i.e.}, the trajectory of the initial condition $pa\in SW_{23}$ reaches the switching plane $SW_{12}$ and not the equilibrium point ${\bf x}^{*}_{eq_2}$. This means that in ${\bf x}$ coordinates, the heteroclinic orbit  from ${\bf x}^{*}_{eq_2}$ to ${\bf x}^{*}_{eq_3}$ does not exist. Similarly, when $pb_z$ is further than $pa_z$ from ${\bf z}^{*(2)}_{eq_2}$, this is, $d(pa_z,{\bf z}^{(2)}_{eq_2})<d(pb_z,{\bf z}^{(2)}_{eq_2})$, for $\gamma$ sufficiently small, the trajectory with initial condition $pb_z$ will eventually reach the set given by \eqref{eq:interseccionp2az} for $t<0$, {\it i.e.}, the trajectory of the initial condition $pb\in SW_{12}$ reaches the switching plane $SW_{23}$ and not the equilibrium point ${\bf x}^{*}_{eq_2}$. This means that in ${\bf x}$ coordinates, the heteroclinic orbit  from ${\bf x}^{*}_{eq_2}$ to ${\bf x}^{*}_{eq_1}$ does not exist.

The next proposition warranty the existence of heteroclinic orbits when $\gamma$ belongs to an interval of real numbers where the case $\gamma-\alpha\neq\alpha$ is considered, such that $\gamma>\alpha$.

\begin{proposition}\label{prop:sixHO}
	The hyperbolic system given by \eqref{eq:affine}, \eqref{eq:Amatrix},\eqref{eq1:vectorB} and \eqref{eq:funcion4atomos} with the switching surfaces given in \eqref{eq:surfaces4atoms} generates six heteroclinic orbits if 
	\[ 
	\frac{\alpha(e^{-a\tau}\cos(b \tau)-1)}{e^{-a\tau}\cos(b \tau)}>\gamma>\alpha(1-e^{-a\tau}\cos(b \tau)),
	\]
	where 
	\[ \tau=\frac{\arctan(b/a)+\pi/2}{b}.\] 
\end{proposition}

{\it Proof:}
To find the values of $\gamma$ for which these heteroclinic orbits exist, let us assume $pa$ is a point of the heteroclinic orbit joining ${\bf x}^{*}_{eq_2}$ and ${\bf x}^{*}_{eq_3}$, {\it i.e.},
\[\lim\limits_{t\to-\infty}\varphi(pa,t)={\bf x}^*_{eq_2}, \ \text{and} \ \lim\limits_{t\to\infty}\varphi(pa,t)={\bf x}^*_{eq_3}.\]
Because $pa\in W^s_{{\bf x}^*_{eq_3}}$ then $\lim\limits_{t\to\infty}\varphi(pa,t)={\bf x}^*_{eq_3}$. For the other part of the heteroclinic orbit, we analyze the system in ${\bf z}^{(2)}$ coordinates, we have $pa_z$, $pb_z$, ${\bf z}^{*(2)}_{eq_2}$, 

and the orbit is given by  ${\bf z}^{(2)}(t)$. 
We assume  that ${\bf z}^{(2)}(0)=pa_z$, so we want that ${\bf z}^{(2)}(t)$  remains in $P_2$ for all $t<0$. 
Thus, we need to find the first maximum in the component $z_3^{(2)}$ of the trajectory whose initial condition is $pa_z$ for $t<0$. According to \eqref{eq:papbz}, the third component of $pa_z$ and $pb_z$ are $\frac{2(\alpha-\gamma)}{3}<0$ and $0<\frac{2\alpha}{3}$, respectively. This maximum give us the intersection point between the trajectory ${\bf z}^{(2)}(t)$ and the axis $z_3^{(2)}$, then we can compare the third components of the trajectory ${\bf z}^{(2)}(t)$ and the point $pb_z$, in terms of $\alpha$ and $\gamma$ to ensure that ${\bf z}^{(2)}(t)$ remains in $P_2$ for all $t<0$.

The trajectory ${\bf z}^{(2)}(t)$ for the initial condition ${\bf z}^{(2)}_0=\left(z^{(2)}_{10},z^{(2)}_{20},z^{(2)}_{30}\right)^T$ is

\begin{equation}\label{eq:z1}
	z^{(2)}_1(t)=z^{(2)}_{10}e^{-c t},
\end{equation}
\begin{equation}\label{eq:z2}
	z^{(2)}_2(t)=z^{(2)}_{20}e^{a t}\cos(b t)-z^{(2)}_{30}e^{a  t}\sin(b t),
\end{equation}
\begin{equation}\label{eq:z3}
	z^{(2)}_3(t)=z^{(2)}_{20}e^{a t}\sin(b t)+z^{(2)}_{30}e^{a  t}\cos(b t),
\end{equation}
this set of equations is analyzed for $t<0$. The same analysis can be done for $0<t$ by using the following set of equations 
\begin{equation}\label{eq:z10}
	z_1^{(2)}(t)=z^{(2)}_{10}e^{c t},
\end{equation}
\begin{equation}\label{eq:z20}
	z_2^{(2)}(t)=z^{(2)}_{20}e^{-a t}\cos(b t)+z^{(2)}_{30}e^{-a  t}\sin(b t),
\end{equation}
\begin{equation}\label{eq:z30}
	z_3^{(2)}(t)=-z^{(2)}_{20}e^{-a t}\sin(b t)+z^{(2)}_{30}e^{-a  t}\cos(b t).
\end{equation}

Since we are looking for the first maximum in $z^{(2)}_3(t)$ for $0<t$.

Then from \eqref{eq:z30} with the initial condition $pa_z$ given in \eqref{eq:papbz}

\begin{equation}\label{eq:paramaximo}
	z^{(2)}_3(t)=\frac{2(\alpha-\gamma)}{3}e^{-a  t}\cos(b t),
\end{equation}
\begin{equation}
	\dot{z}^{(2)}_3(t)=-\frac{2(\alpha-\gamma)}{3}e^{-a  t}\left(b\sin(bt)+a\cos(b t)\right),
\end{equation}
\begin{equation}
	\dot{z}^{(2)}_3(t)=-\frac{2(\alpha-\gamma)}{3}e^{-a  t}\left(\sqrt{a^2+b^2}\cos(b t-\arctan(b/a))\right),
\end{equation}
to find the maximum we equate to zero
\begin{equation}
	0=-\frac{2(\alpha-\gamma)}{3}e^{-a  t}\left(\sqrt{a^2+b^2}\cos(b t-\arctan(b/a))\right),
\end{equation}
thus, it turns out that
\begin{equation}
	b t-\arctan(b/a)=\pi/2+n\pi, \mbox{ with }n\in\mathbb{Z},
\end{equation}
\begin{equation}\label{eq:tiempo}
	t=\frac{\arctan(b/a)}{b}+\frac{\pi}{2b}+\frac{n\pi}{b}, \mbox{ with }n\in\mathbb{Z}.
\end{equation}
We will call $t_{max}$ the time for the first maximum. 
Thus it follows that
\begin{equation}
	t_{max}=\frac{\arctan(b/a)+\pi/2}{b},
\end{equation}

then from \eqref{eq:paramaximo} 
\begin{equation}
	z^{(2)}_3(t_{max})=\frac{2(\alpha-\gamma)}{3}e^{-a  t_{max}}\cos(b t_{max}),
\end{equation}
this maximum $z_3^{(2)}$ must be part of $P_2$, since $pb_z$ belongs to $P_1$ it follows from \eqref{eq:papbz} that
\begin{equation}
	\frac{2\alpha}{3}>\frac{2(\alpha-\gamma)}{3}e^{-at_{max}}\cos(b t_{max}),
\end{equation}
\begin{equation}
	\frac{\alpha(e^{-at_{max}}\cos(b t_{max})-1)}{e^{-at_{max}}\cos(b t_{max})}>\gamma.
\end{equation}

Now, let us assume $pb$ is a point of the heteroclinic orbit joining ${\bf x}^{*}_{eq_2}$ and ${\bf x}^{*}_{eq_1}$, {\it i.e.},
\[\lim\limits_{t\to-\infty}\varphi(pb,t)={\bf x}^*_{eq_2}, \ \text{and} \ \lim\limits_{t\to\infty}\varphi(pb,t)={\bf x}^*_{eq_1}.\]
Because $pb\in W^s_{{\bf x}^*_{eq_1}}$ then $\lim\limits_{t\to\infty}\varphi(pb,t)={\bf x}^*_{eq_1}$. 

Following the same procedure described above but looking for a minimum, due to the third component of $pb_z$ is $0<\frac{2\alpha}{3}$.  it is found that
\begin{equation}
	t_{min}=\frac{\arctan(b/a)+\pi/2}{b},
\end{equation}
then from \eqref{eq:z30} and the point $pb_z$ given in \eqref{eq:papbz}
\begin{equation}
	z^{(2)}_3(t_{min})=\frac{2\alpha}{3}e^{-at_{min}}\cos(b t_{min}),
\end{equation}
this minimum $z_3^{(2)}$ must be part of $P_2$, since $pa_z$ belongs to $P_3$ it follows from \eqref{eq:papbz} that
\begin{equation}
	\frac{2(\alpha-\gamma)}{3}<\frac{2\alpha}{3}e^{-at_{min}}\cos(b t_{min}),
\end{equation}
\begin{equation}
	\gamma>\alpha(1-e^{-at_{min}}\cos(b t_{min})).
\end{equation}
Then defining $\tau=t_{max}=t_{min}$
\begin{equation}
	\frac{\alpha(e^{-a\tau}\cos(b \tau)-1)}{e^{-a\tau}\cos(b \tau)}>\gamma>\alpha(1-e^{-a\tau}\cos(b \tau)).
\end{equation}
The same conclusion apply to the point ${\bf x}^*_{eq_3}$ due to the symmetry of the system. Finally, the heteroclinic orbit from ${\bf x}^*_{eq_1}$ to ${\bf x}^*_{eq_2}$ and  the one from ${\bf x}^*_{eq_4}$ to ${\bf x}^*_{eq_3}$ are always present in the system as there are no more switching surfaces.   $\qed$

\begin{figure}[h!]
	\centering
	\subfloat[]{\includegraphics[width=0.45\columnwidth]{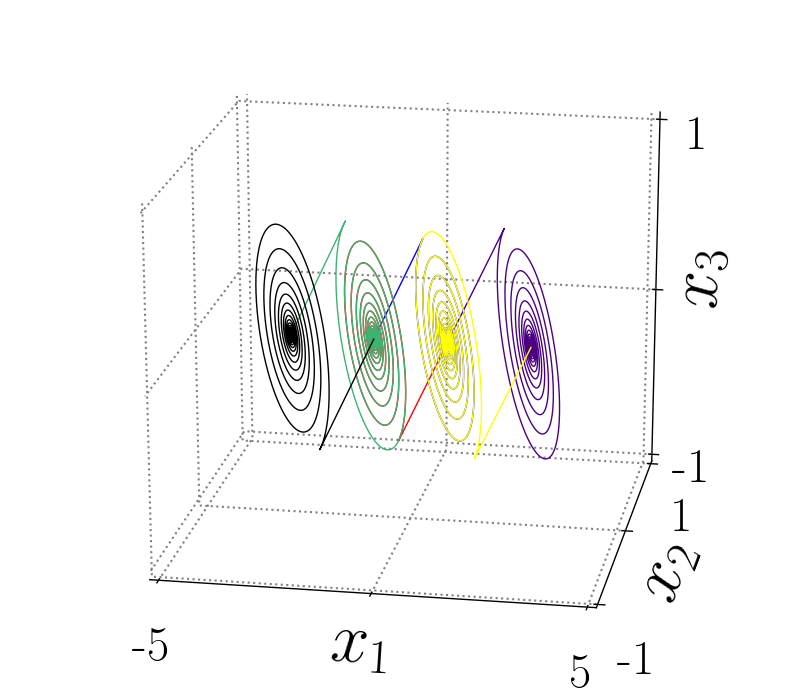}\label{fig:mindelta}}
	\subfloat[]{\includegraphics[width=0.45\columnwidth]{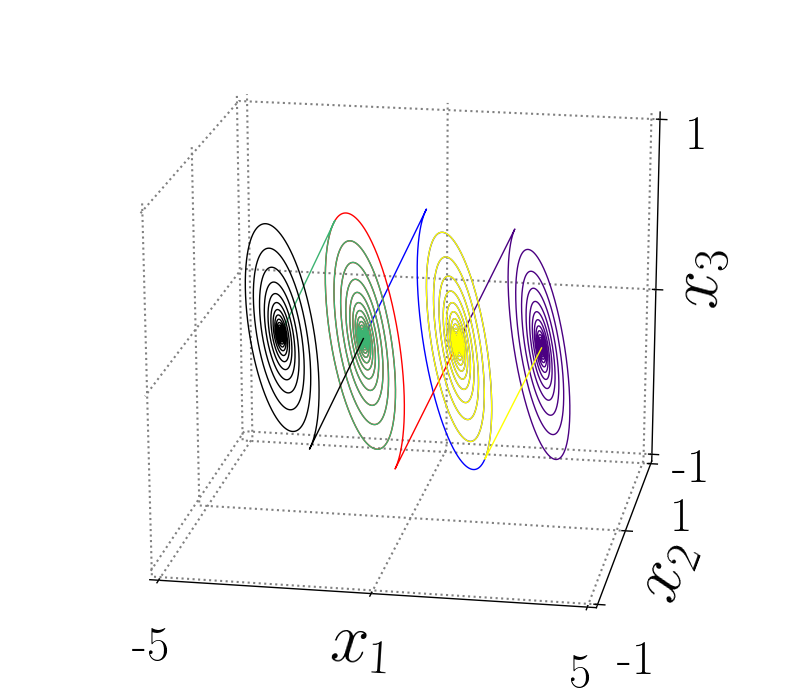}\label{fig:maxdelta}}\\
	\subfloat[]{\includegraphics[width=0.45\columnwidth]{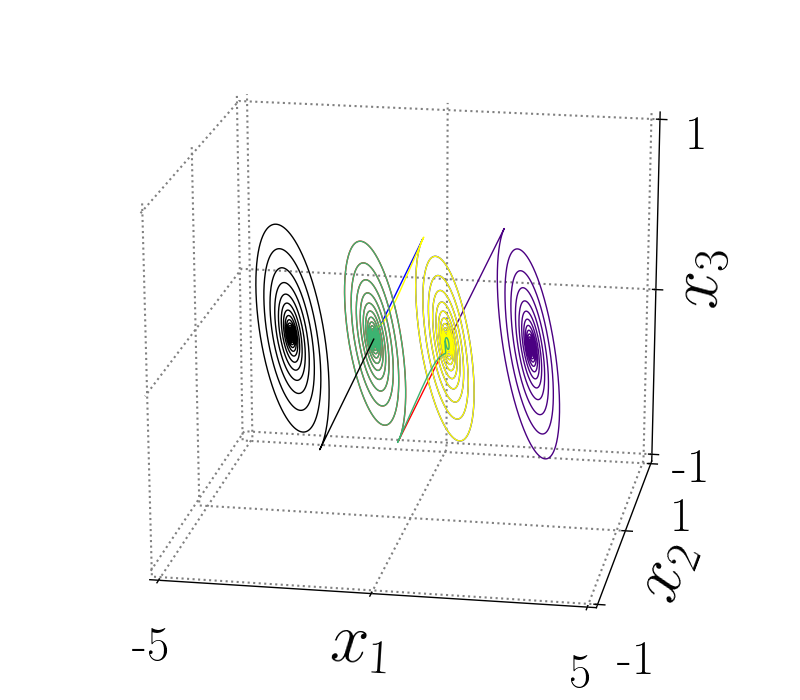}\label{fig:min}}
	\subfloat[]{\includegraphics[width=0.45\columnwidth]{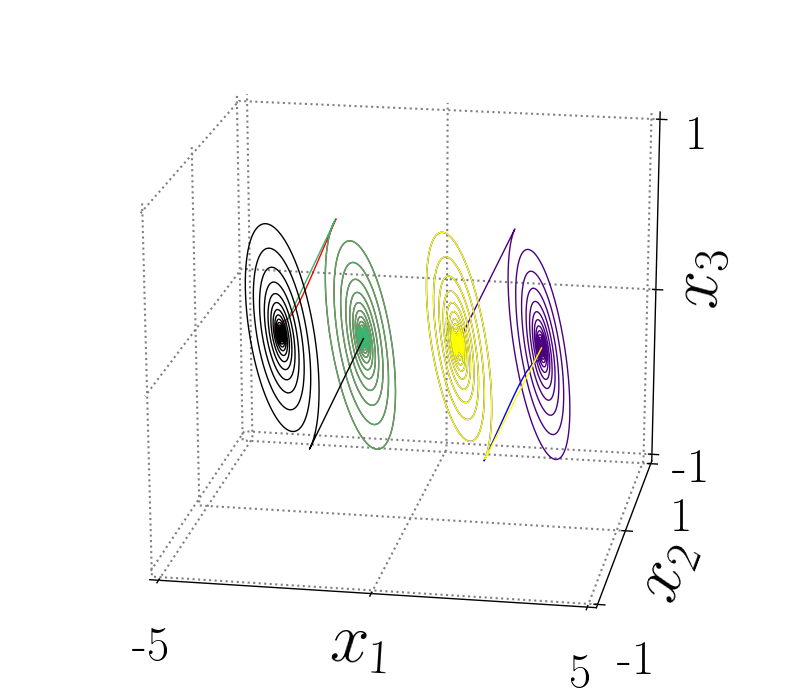}\label{fig:max}}
	\caption{\label{fig:casesorbits}Heteroclinic orbits of the system given by \eqref{eq:affine}, \eqref{eq:Amatrix}, \eqref{eq1:vectorB}, \eqref{eq:funcion4atomos} and \eqref{eq:surfaces4atoms} for the parameters $a=0.2$, $b=5$, $c=-3$, $\alpha=1$ and different values of $\gamma$. There are six heteroclinic orbits for: (a)   $\gamma_1=\alpha(1-e^{-a\tau}\cos(b \tau))+.00001$, and (b)  $\gamma_2=\frac{\alpha(e^{-a\tau}\cos(b \tau)-1)}{e^{-a\tau}\cos(b \tau)}-.00001$. Four heteroclinic orbits for: (c)  $\gamma_L=\alpha(1-e^{-a\tau}\cos(b \tau))$, and (d) $\gamma_U=\frac{\alpha(e^{-a\tau}\cos(b \tau)-1)}{e^{-a\tau}\cos(b \tau)}$.  }
\end{figure}

To illustrate the effect of the parameters $\gamma$, $\alpha$, $a$, and $b$ on the existence of heteroclinic orbits of the system given by \eqref{eq:affine}, \eqref{eq:Amatrix}, \eqref{eq1:vectorB}, \eqref{eq:funcion4atomos} and \eqref{eq:surfaces4atoms}, we use the proposition \eqref{prop:sixHO} to determine the open interval of real values for $\gamma$  given by 
$$\Gamma=\left(\alpha(1-e^{-a\tau}\cos(b \tau)),\frac{\alpha(e^{-a\tau}\cos(b \tau)-1)}{e^{-a\tau}\cos(b \tau)}\right),$$
with $\tau=\frac{\arctan(b/a)+\pi/2}{b}$.
So, six initial conditions were calculated as in \eqref{eq:x01} and \eqref{eq:x02} with $k=50$ for the parameters $a=0.2$, $b=5$, $c=-3$, and $\alpha=1$. Four cases of different values of $\gamma$ are analyzed. The first two correspond to $\gamma_{1,2}\in \Gamma$ and the last two correspond to $\gamma_{L,U}\notin \Gamma$
:
\begin{itemize}
	\item[1.-] For this case $\gamma_{1}\in \Gamma$, with $\gamma_1=\alpha(1-e^{-a\tau}\cos(b \tau))+.00001$, so there exist {\it six heteroclinic orbits} as shown in the Figure~\ref{fig:mindelta}.
	\item[2.-] For $\gamma_2=\left(\frac{\alpha(e^{-a\tau}\cos(b \tau)-1)}{e^{-a\tau}\cos(b \tau)}-.00001 \right) \in \Gamma$, in this case there exist also {\it six heteroclinic orbits} as shown in the Figure~\ref{fig:maxdelta}.
	\item[ 3.-] In this case $\gamma_L=\alpha\left( 1-e^{-a\tau}\cos(b \tau) \right) \notin \Gamma$,  then there exist {\it four heteroclinic orbits} as shown in Figure~\ref{fig:min}. The green orbit stating close to ${\bf x}^*_{eq_2}$ cannot reach ${\bf x}^*_{eq_1}$ and goes to $P_3$. In the same way, the yellow orbit starting close to ${\bf x}^*_{eq_3}$ cannot reach ${\bf x}^*_{eq_4}$ and goes to $P_2$. Then there is no heteroclinic orbits from  ${\bf x}^*_{eq_2}$ to ${\bf x}^*_{eq_1}$ and from   ${\bf x}^*_{eq_3}$ to ${\bf x}^*_{eq_4}$.
	\item[ 4.-] For $\gamma_U=\frac{\alpha(e^{-a\tau}\cos(b \tau)-1)}{e^{-a\tau}\cos(b \tau)} \notin \Gamma$, there exist also {\it four heteroclinic orbits} as shown in the Figure~\ref{fig:max}. The red orbit stating close to ${\bf x}^*_{eq_2}$ cannot reach ${\bf x}^*_{eq_3}$ and goes to $P_1$. In the same way, the blue orbit starting close to ${\bf x}^*_{eq_3}$ cannot reach ${\bf x}^*_{eq_2}$ and goes to $P_4$. Then there is no heteroclinic orbit from  ${\bf x}^*_{eq_2}$ to ${\bf x}^*_{eq_3}$, nor vice versa.
\end{itemize}

\begin{figure}[!ht]
	\centering
	\subfloat[]{\includegraphics[width=0.3\columnwidth]{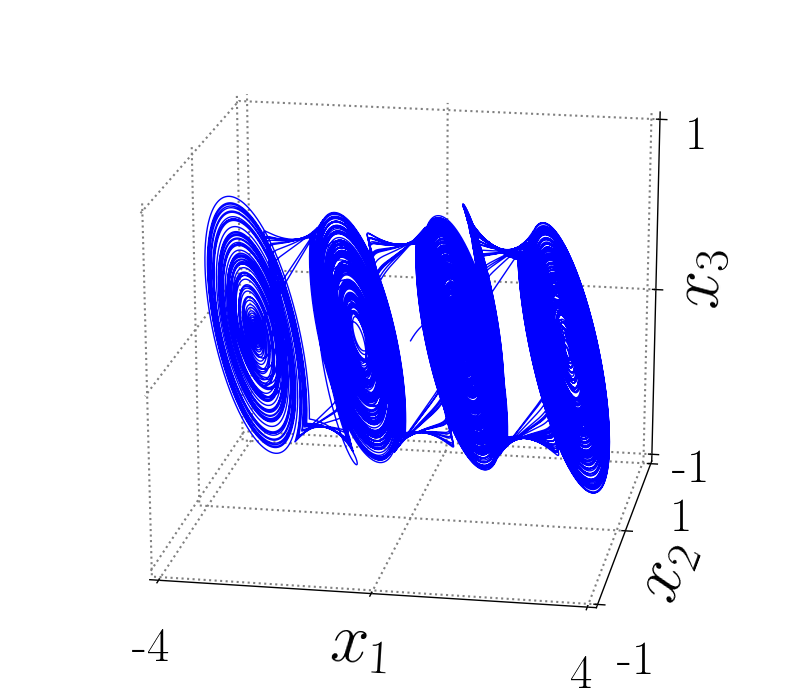}\label{fig:multiscroll}}
	\subfloat[]{\includegraphics[width=0.3\columnwidth]{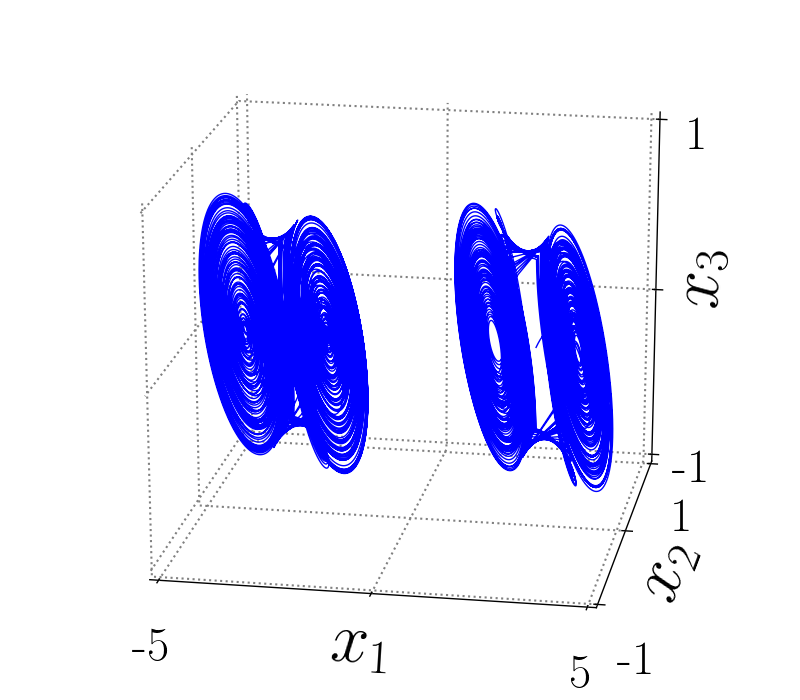}\label{fig:separated}}
	\subfloat[]{\includegraphics[width=0.3\columnwidth]{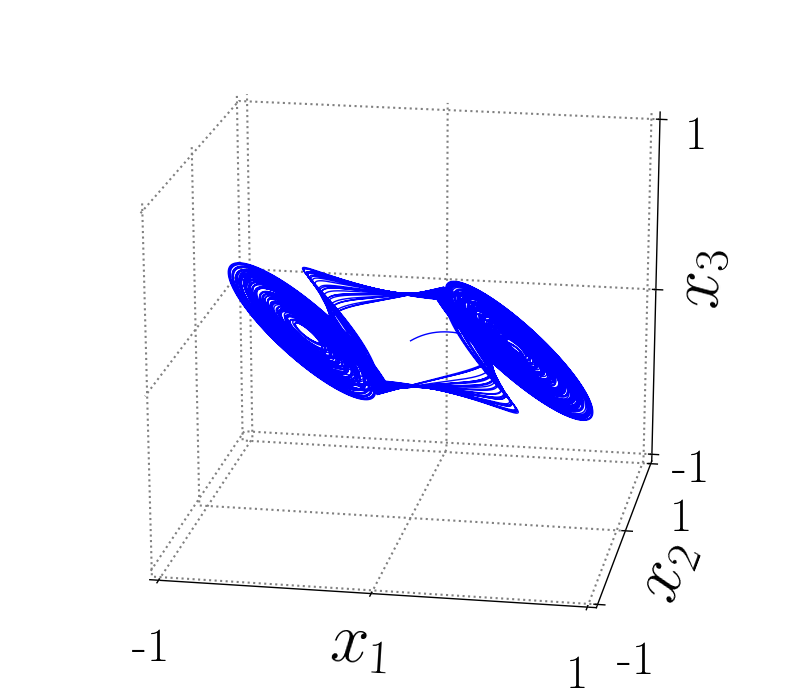}\label{fig:central}}\\
	\subfloat[]{\includegraphics[width=0.3\columnwidth]{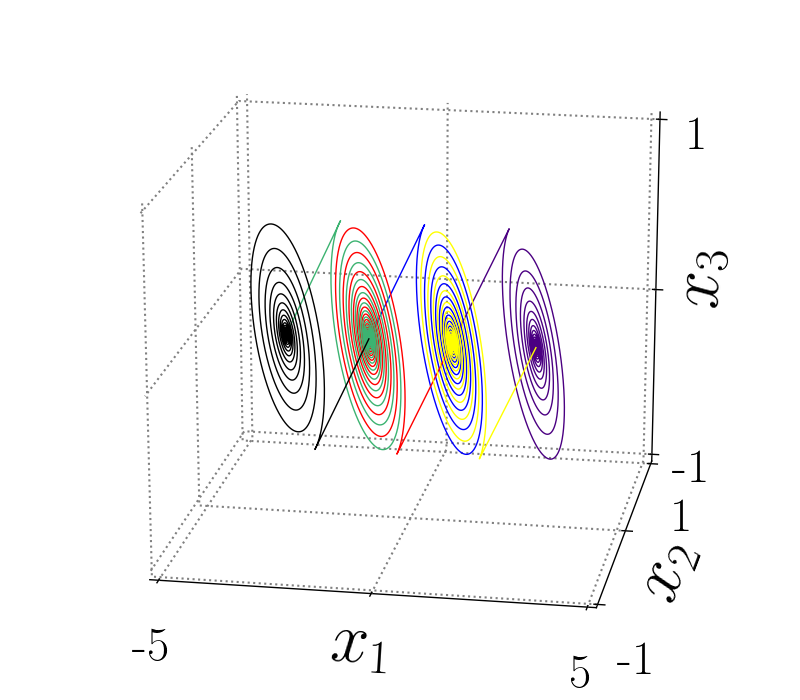}\label{fig:multiscrollorbits}}
	\subfloat[]{\includegraphics[width=0.3\columnwidth]{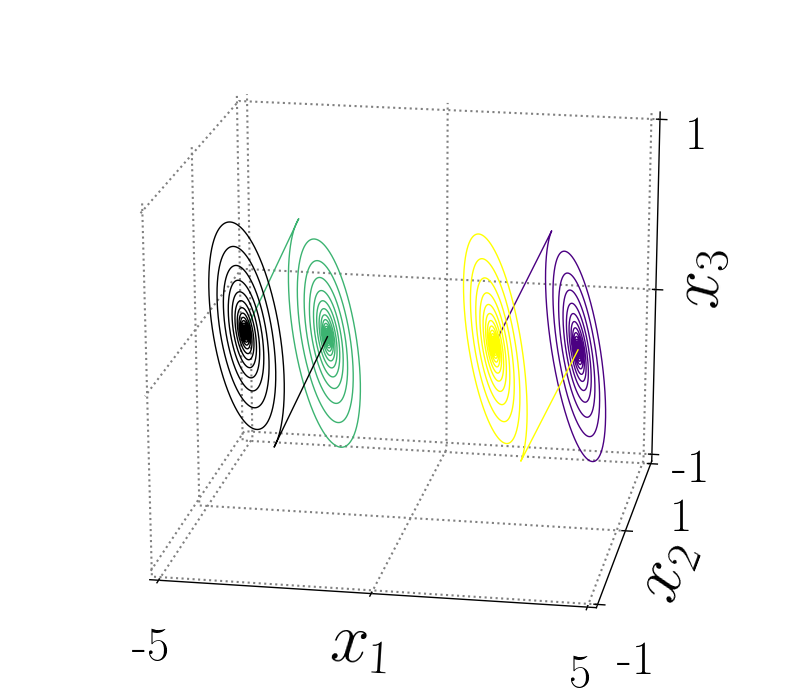}\label{fig:separatedorbits}}
	\subfloat[]{\includegraphics[width=0.3\columnwidth]{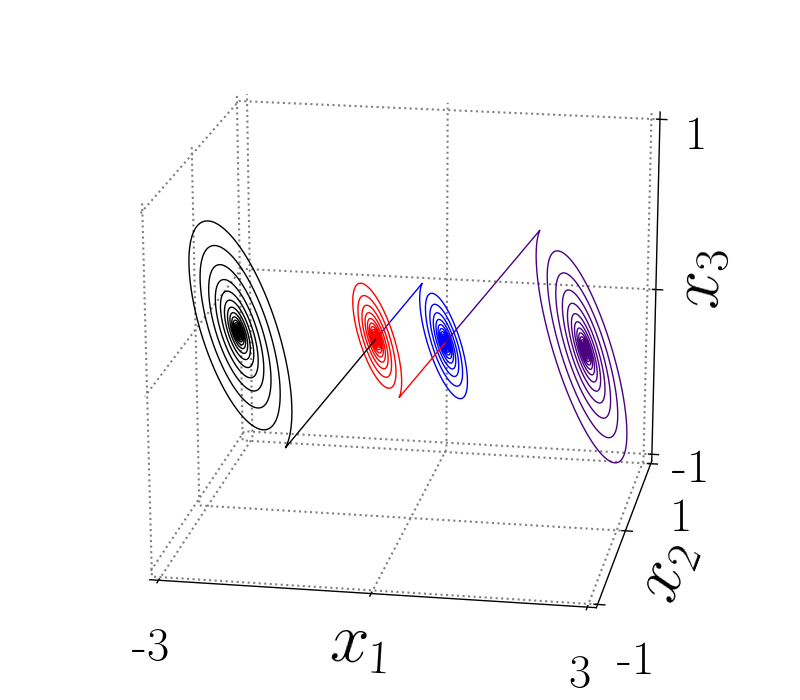}\label{fig:centralorbits}}
	\caption{\label{fig:cases}Attractors and heteroclinic orbits of the system given by \eqref{eq:affine}, \eqref{eq:Amatrix}, \eqref{eq1:vectorB}, \eqref{eq:funcion4atomos} and \eqref{eq:surfaces4atoms} for $\alpha=1$ and different values of $\gamma$. For $\gamma=2$ the system exhibits in (a) a quad-scroll attractor, and in (d) six heteroclinic orbits. For $\gamma=3$ the system exhibits in (b) two double scroll attractors, and in (e) four heteroclinic orbits. For $\gamma=1.5$ the system exhibits in (c) a double-scroll attractor, and in (f) four heteroclinic orbits.}
\end{figure}

The open interval $\Gamma$ is given as  $\Gamma=(\gamma_L,\gamma_U)$, where
\begin{equation}
	\gamma_L=\alpha(1-e^{-a\tau}\cos(b \tau))\approx1.8826170015164836,
\end{equation}
\begin{equation}
	\gamma_U=\frac{\alpha(e^{-a\tau}\cos(b \tau)-1)}{e^{-a\tau}\cos(b \tau)}\approx2.1329942639693464.
\end{equation}

The four cases mentioned generate three types of systems determined by $\gamma$ and $\Gamma$. For instance, for $\gamma=2\in\Gamma$ and $\alpha=1$ corresponds to the above first and second cases. Then the system presents six heteroclinic orbits which comprise three heteroclinic loops between equilibria: ${\bf x}^*_{eq_1}$ and ${\bf x}^*_{eq_2}$; ${\bf x}^*_{eq_2}$ and ${\bf x}^*_{eq_3}$; ${\bf x}^*_{eq_3}$ and ${\bf x}^*_{eq_4}$. For $\gamma=1.5<\gamma_L$, then  $\gamma\notin\Gamma$, and this case corresponds to the above third case. So there are four heteroclinic orbits and two of them comprise a heteroclinic loop between equilibria ${\bf x}^*_{eq_2}$ and ${\bf x}^*_{eq_3}$.  For $\gamma_U<\gamma=3$, then  $\gamma\notin\Gamma$, and this case corresponds to the above fourth case. So there are four heteroclinic orbits which comprise two heteroclinic loops, but now between equilibria: ${\bf x}^*_{eq_1}$ and ${\bf x}^*_{eq_2}$;  ${\bf x}^*_{eq_3}$ and ${\bf x}^*_{eq_4}$. 
The above three cases generate self-excited attractors as shown below: 
\begin{itemize}
	\item[1.-] For $\gamma=2\in \Gamma$, the system presents a self-excited attractor with four scrolls which is shown in the Figure~\ref{fig:multiscroll} and its corresponding  three heteroclinic loops are shown in the  Figure~\ref{fig:multiscrollorbits}.
	According to \cite{Campos-Canton10} a scroll attractor with three or more scrolls is considered as a multiscroll attractor, thus the attractor shown in the Figure~\ref{fig:multiscroll} is a multiscroll attractor. The scrolls are generated around each equilibrium point of the system ${\bf x}^*_{eq_i}$, with $i=1,2,3,4$.
	\item[2.-] For $\gamma=3$, $\gamma>\gamma_U$, then $\gamma\notin \Gamma$. The system presents bistability, the two double-scroll self-excited attractors are shown in the Figure~\ref{fig:separated}. In this case two heteroclinic orbits are lost, the system exhibits four heteroclinic orbits, {\it i.e.}, two heteroclinic loops, as shown in Figure~\ref{fig:separatedorbits}. One  double-scroll self-excited attractor oscillates around equilibria  ${\bf x}^*_{eq_1}$ and ${\bf x}^*_{eq_2}$, while the other self-excited attractor oscillates around equilibria  ${\bf x}^*_{eq_3}$ and ${\bf x}^*_{eq_4}$. The basin of attraction of each self-excited attractor is surrounded both attractors. 
	\item[3.-] For $\gamma=1.5$, $\gamma<\gamma_L$, then $\gamma\notin \Gamma$.  The system presents only one double-scroll self-excited attractor shown in the Figure~\ref{fig:central}. In this case two heteroclinic orbits are also lost but only a heteroclinic loop is exhibited. The heteroclinic orbits are shown in the Figure~\ref{fig:centralorbits}.  The double-scroll self-excited attractor oscillates around equilibria  ${\bf x}^*_{eq_2}$ and ${\bf x}^*_{eq_3}$,
\end{itemize}

Based on the results of widening of the basins of attraction of a
multistable switching dynamical system with the location of symmetric equilibria
reported in \cite{Ontanon17}, we could ponder in the possible existence of a hidden attractor for the case  $\gamma>\gamma_U$, because there are oscillations surround the two self-exited attractors as a hidden attractor exists. However, the simulations of these systems let us know that hidden attractors are not present. For example, if the $\gamma$ value is increased  then also the distance between the two self-exited attractors increases. This provokes that some initial conditions in the basins of attraction of both attractors generate transitory oscillations resembling a double scroll attractor, however, after some time these transitory oscillations converge to one of the double-scroll self-excited attractors.

\begin{figure}[!ht]
	\centering
	\subfloat[]{\includegraphics[width=0.3\columnwidth]{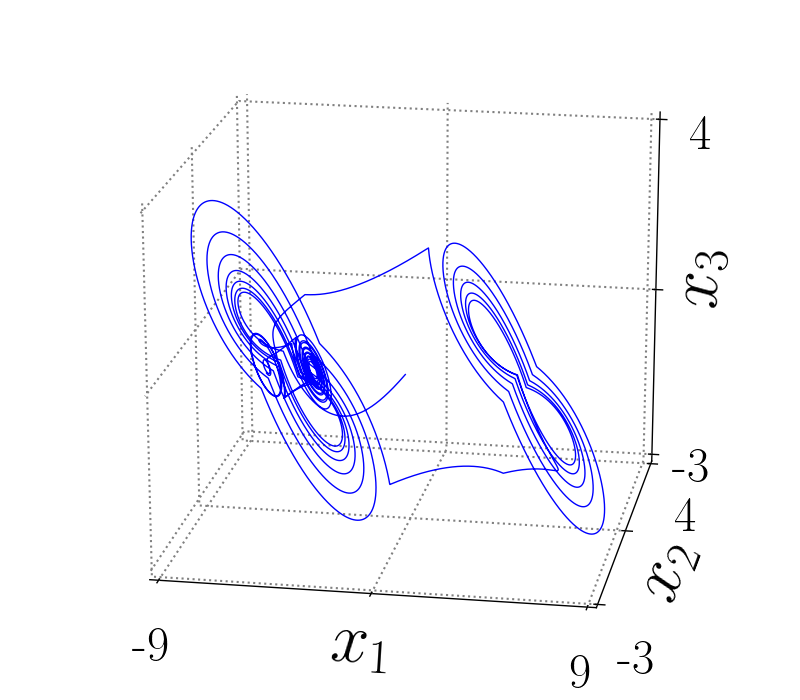}\label{fig:transitory5}}
	\subfloat[]{\includegraphics[width=0.3\columnwidth]{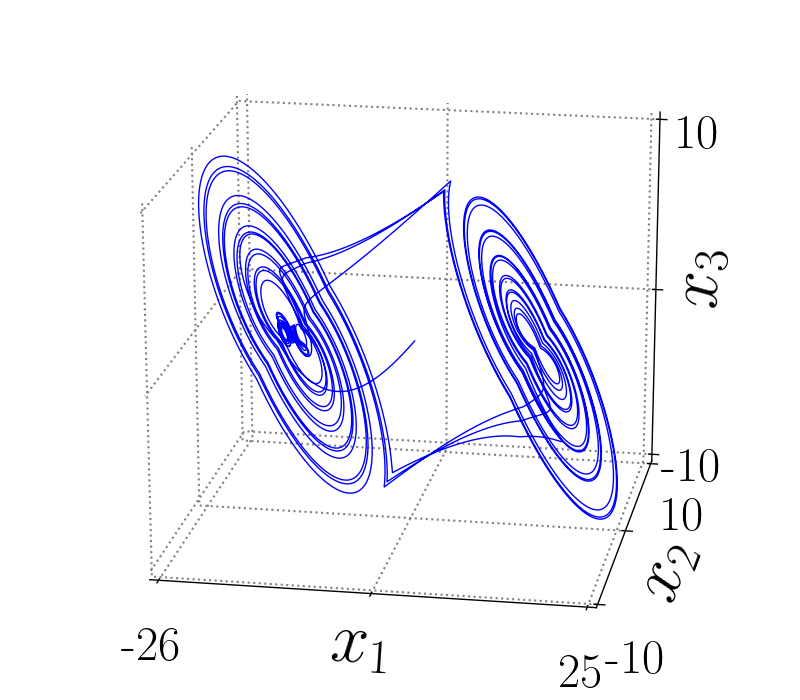}\label{fig:transitory15}}
	\subfloat[]{\includegraphics[width=0.3\columnwidth]{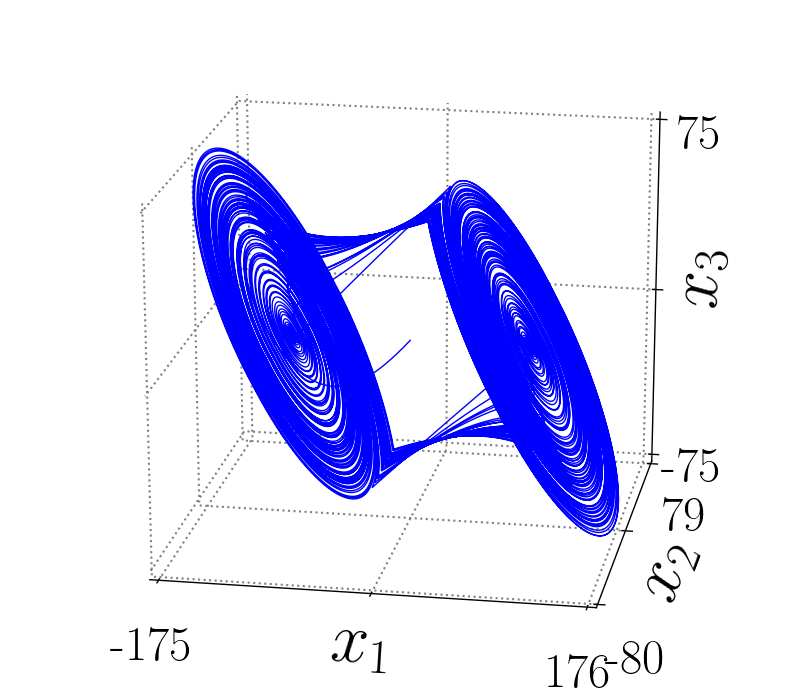}\label{fig:transitory100}}\\
	\subfloat[]{\includegraphics[width=0.3\columnwidth]{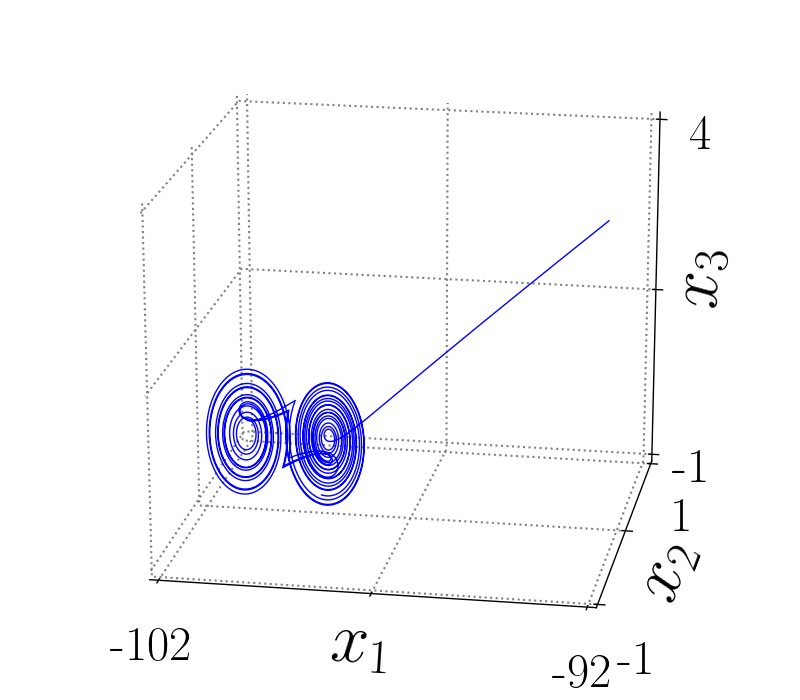}\label{fig:transitory100b}}
	\subfloat[]{\includegraphics[width=0.3\columnwidth]{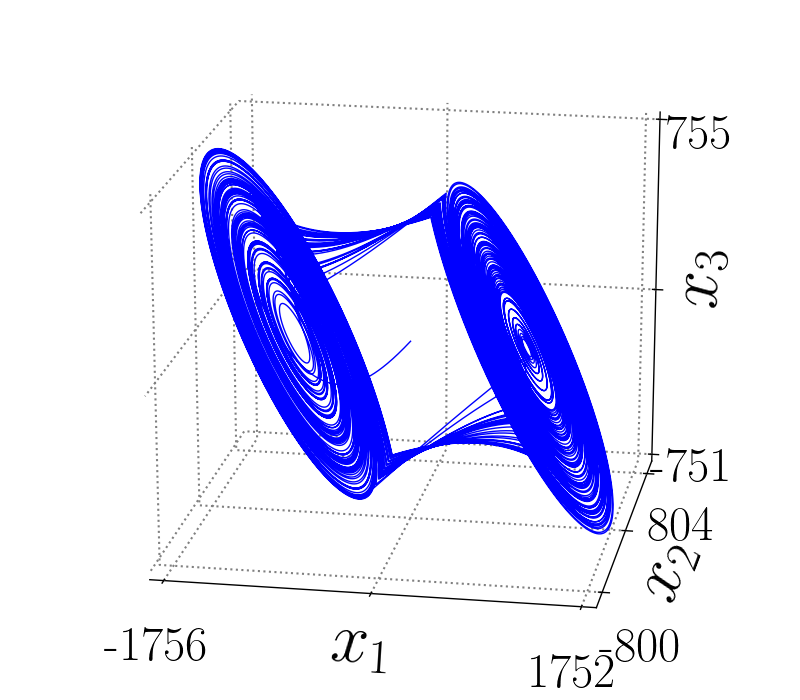}\label{fig:transitory1000}}
	\subfloat[]{\includegraphics[width=0.3\columnwidth]{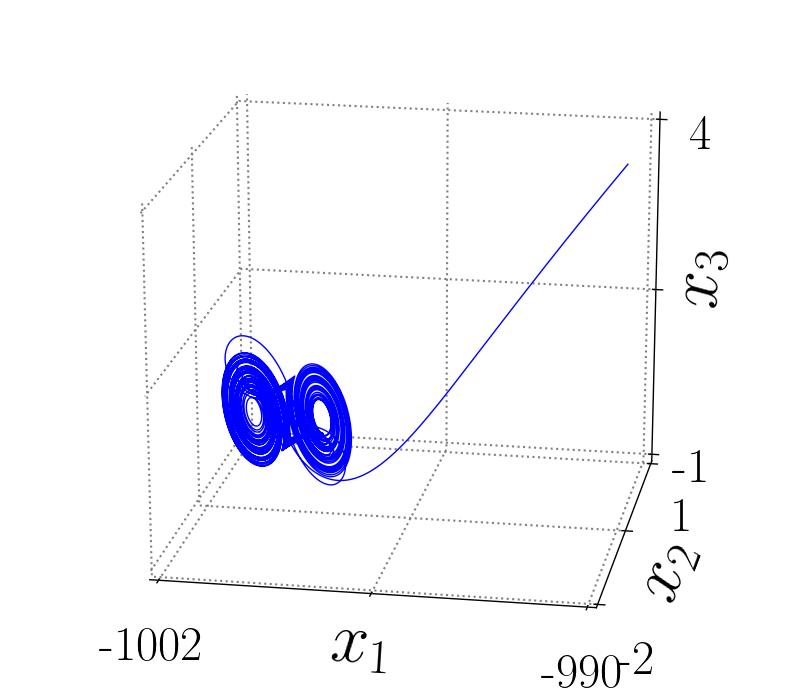}\label{fig:transitory1000b}}
	\caption{\label{fig:transitory} Trajectory of the system given by \eqref{eq:affine}, \eqref{eq:Amatrix}, \eqref{eq1:vectorB}, \eqref{eq:funcion4atomos} and \eqref{eq:surfaces4atoms} for the initial condition $x=(0,0,0)^T$, $a=0.2$, $b=5$, $c=-7$, $\alpha=1$ and different values of $\gamma$: (a) $\gamma=5$, the transitory oscillation of double-scroll exhibited and after some time converge to one of the double-scroll self-excited attractors, $t\in[0,40]$ $a.u.$; (b)$\gamma=15$, the transitory oscillation of double-scroll exhibited and after some time converge to one of the double-scroll self-excited attractors, $t\in[0,60]$ $a.u.$; (c) $\gamma=100$, the transitory oscillation of double scroll exhibited for $t\in[0,300]$ $a.u.$; (d) $\gamma=100$, double-scroll self-excited attractor for $t=[356.6,400]$ $a.u.$; (e) $\gamma=1000$, the transitory oscillation of double scroll exhibited for $t\in[0,300]$ $a.u.$; (f)  $\gamma=1000$, the double-scroll self-excited attractor for $t=[3091,3200]$ $a.u.$.    }
\end{figure}

We analyze the trajectory for the initial condition ${\bf x}_0=(0,0,0)^T$ and different values of $\gamma$ fulfilling  $\gamma_U<\gamma$. The first case is $\gamma=5$ and  $t\in[0,40]$ in arbitrary units ($a.u.$). The equilibria are at:
\begin{equation}
	{\bf x}^*_{eq_1}=\begin{bmatrix}
		-6\\0\\0
	\end{bmatrix},\ 
	{\bf x}^*_{eq_2}=\begin{bmatrix}
		-4\\0\\0
	\end{bmatrix},\
	{\bf x}^*_{eq_3}=\begin{bmatrix}
		4\\0\\0
	\end{bmatrix},\
	{\bf x}^*_{eq_4}=\begin{bmatrix}
		6\\0\\0
	\end{bmatrix}.
\end{equation}
Figure~\ref{fig:transitory5} shows the trajectory which consists of the transitory behavior resembles a double scroll attractor and after a short time reaches a double-scroll self-excited attractor around equilibria ${\bf x}^*_{eq_1}$ and ${\bf x}^*_{eq_2}$.

Increasing the value of $\gamma$ to $15$, the equilibria are located at:
\begin{equation}
	{\bf x}^*_{eq_1}=\begin{bmatrix}
		-16\\0\\0
	\end{bmatrix},\ 
	{\bf x}^*_{eq_2}=\begin{bmatrix}
		-14\\0\\0
	\end{bmatrix},\
	{\bf x}^*_{eq_3}=\begin{bmatrix}
		14\\0\\0
	\end{bmatrix},\
	{\bf x}^*_{eq_4}=\begin{bmatrix}
		16\\0\\0
	\end{bmatrix}.
\end{equation}
And the transitory time to reaches the self-exited attractor  is increased.
In the Figure~\ref{fig:transitory15} the trajectory is shown for $t\in [0,60]$ $a.u.$.
Now, for $\gamma =100$ the equilibria are located at:
\begin{equation}
	{\bf x}^*_{eq_1}=\begin{bmatrix}
		-101\\0\\0
	\end{bmatrix},\ 
	{\bf x}^*_{eq_2}=\begin{bmatrix}
		-99\\0\\0
	\end{bmatrix},\
	{\bf x}^*_{eq_3}=\begin{bmatrix}
		99\\0\\0
	\end{bmatrix},\
	{\bf x}^*_{eq_4}=\begin{bmatrix}
		101\\0\\0
	\end{bmatrix}.
\end{equation}
The transitory time last longer for the same initial condition. Figure~\ref{fig:transitory100} shows the transitory oscillations of the trajectory for $t\in [0,300]$ $a.u.$. After a long time, the trajectory reaches a double-scroll self-excited attractor around equilibria ${\bf x}^*_{eq_1}$ and ${\bf x}^*_{eq_2}$, see Figure~\ref{fig:transitory100b}  for $t\in[356.6,400]$ $a.u.$.
Continuing to increase the value to $\gamma=1000$, then this sets the equilibria at:
\begin{equation}
	{\bf x}^*_{eq_1}=\begin{bmatrix}
		-1001\\0\\0
	\end{bmatrix},\ 
	{\bf x}^*_{eq_2}=\begin{bmatrix}
		-999\\0\\0
	\end{bmatrix},\
	{\bf x}^*_{eq_3}=\begin{bmatrix}
		999\\0\\0
	\end{bmatrix},\
	{\bf x}^*_{eq_4}=\begin{bmatrix}
		1001\\0\\0
	\end{bmatrix}.
\end{equation}
Figures~\ref{fig:transitory1000} shows the transitory oscillation of the trajectory when $t\in [0,300]$ $a.u.$, again transitory time increases and after this long time, the trajectory again reaches a double-scroll self-excited attractor around equilibria ${\bf x}^*_{eq_1}$ and ${\bf x}^*_{eq_2}$, see \ref{fig:transitory1000b} for $t\in[3091,3200]$ $a.u.$. 

In brief, for $\gamma=5$ it took the trajectory around $35$ $a.u.$ to converges to a self excited attractor, for $\gamma=15$ around $50$ $a.u.$ to converge, for $\gamma=100$ around $350$ $a.u.$ and for $\gamma=1000$ around $3090$ $a.u.$. Thus transitory time seems to increase when $\gamma$ increases. In all the cases the trajectories reach a self-excited attractor. 

There are two interests, the former is to understand how the transient oscillation of the trajectory that resembles a double scroll attractor is directed to a self-excited attractor. The latter is to devise a mechanism to block the flow towards a self-excited attractor in order to generate a hidden attractor.

%%%%%%%%%%%%%
%%Section 4%%
%%%%%%%%%%%%%
\section{Route to a self-excited attractor}\label{Sec:SEA}
\begin{figure}[!ht]
	%\centering
	\subfloat[]{\includegraphics[width=0.4\columnwidth]{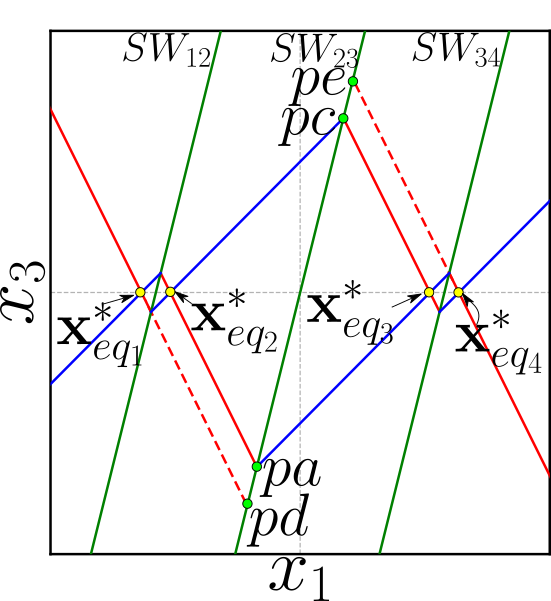}\label{fig:diagrama1}}		 \hspace{1cm}
	\subfloat[]{\includegraphics[width=0.45\columnwidth]{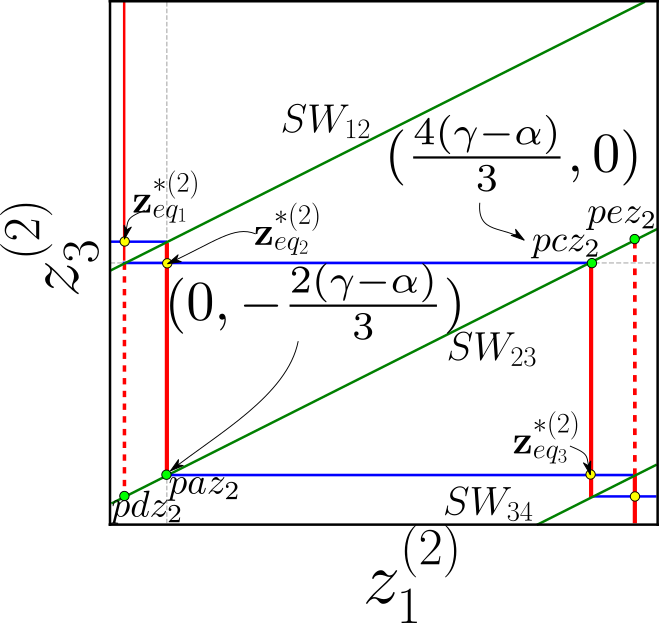}\label{fig:diagramaz2}}
	\hfill
	\caption{\label{fig:diagramas} Projection of the manifolds on (a) $(x_1-x_3)$ plane and (b) $(z^{(2)}_1-z^{(2)}_3)$ plane. The stable and unstable manifolds are marked with blue and red solid lines, respectively, the switching surfaces with green lines.}
\end{figure}

The trajectories with initial conditions in the points $pa, pc \in SW_{23}$ converge to equilibria ${\bf x}_{eq_3}$ and ${\bf x}_{eq_2}$, respectively. So, the transient oscillation of the trajectory that resembles a double scroll attractor is interfered when the trajectory reaches neighborhoods $N(pa)\subset SW_{23}$ and $N(pc)\subset SW_{23}$ around $pa$ or $pc$, respectively, because each trajectory with initial condition in $N(pa)$ or $N(pc)$ is led to one of the self-excited attractors  $A_{self1}$ or $A_{self2}$, respectively. 
Therefore, there are two regions $R_1, R_2 \subset SW_{23}$ inside that can interfere the transient oscillation of the trajectory when $N(pa)\subset R_1$ and $N(pc)\subset R_2$. Thus, these regions contain the points $pa\in R_1$ and $pc\in R_2$. So, the aim of this Section \ref{Sec:SEA} is to understand the location of these regions $R_1, R_2 \subset SW_{23}$ and visualize the route to a self-excited attractor when $N(pa)\cap R_1\neq \emptyset$ and $N(pc)\cap R_2\neq \emptyset$.

In order to understand the transitory behavior, the vector field of $P_2$ and $P_3$ in $SW_{23}$ is analyzed. Then, we need to find a region $R_1 \subset SW_{23}$, such that any trajectory $\varphi(x_0)$, with $x_0\in R_1$, will eventually
go to the self excited attractor $A_{self2}$. The regions $R_1$ and $R_2$ are symmetric with respect to the origin, for instance, the point $pa$ and its symmetric point $pc=-pa$:

\begin{equation}\label{eq:paypb}
	pa=\begin{pmatrix}
		\displaystyle -\frac{\gamma-\alpha}{3}\\0\\\displaystyle-\frac{2(\gamma-\alpha)}{3}\\
	\end{pmatrix},\
	pc=\begin{pmatrix}
		\displaystyle \frac{\gamma-\alpha}{3}\\0\\\displaystyle\frac{2(\gamma-\alpha)}{3}\\
	\end{pmatrix}.
\end{equation}

\begin{figure}[!ht]
	\centering
	\includegraphics[width=\columnwidth]{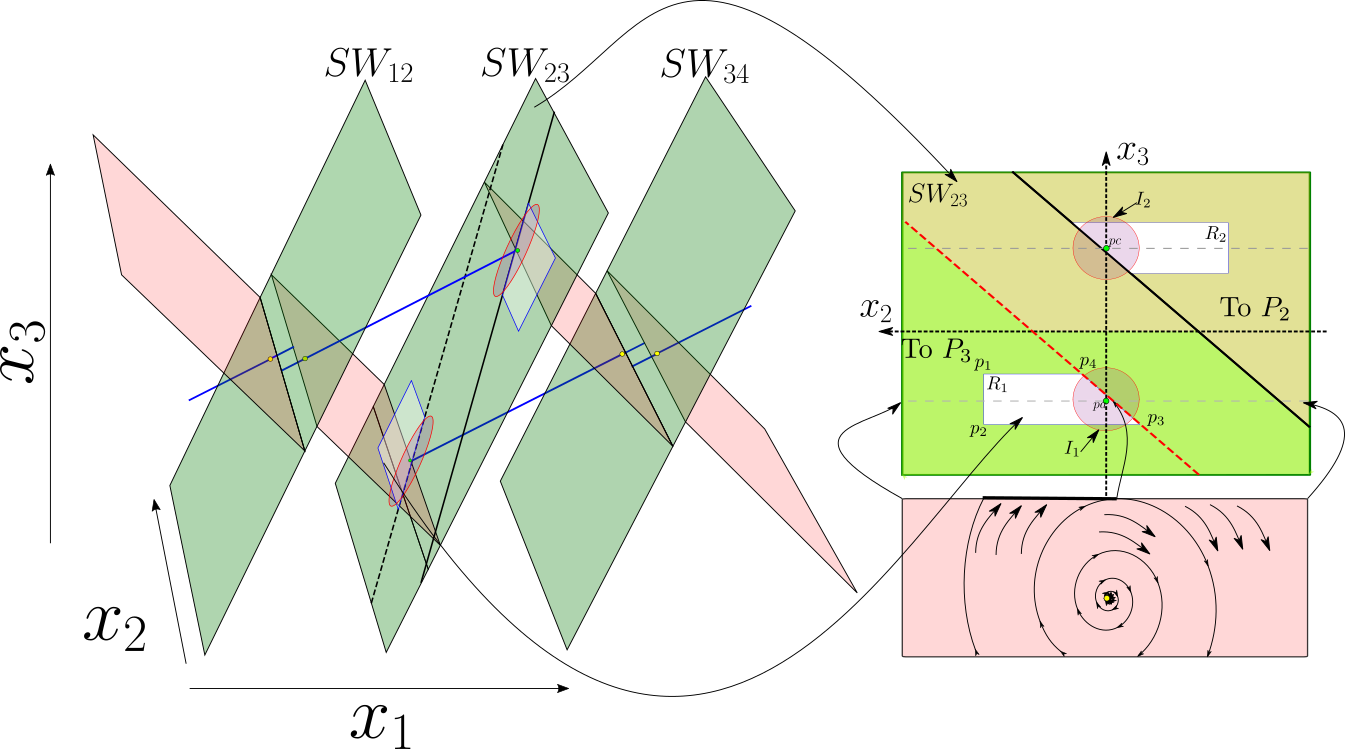}
	\caption{\label{fig:diagrama3d}Illustration of the local manifold of the system. Stable manifolds are in blue, unstable manifold in red and switching surfaces in green. The points where the vector field of $P_2$ is tangent to $SW_{23}$ are in a red dotted line while the points where the vector field of $P_3$ are tangent to $SW_{23}$ are in a continuous black line. The symmetric regions $R_1$ and $R_2$ are in white.}
\end{figure}

To start the analysis, let us find the points in $SW_{23}$ where the vector field of $P_2$ and $P_3$ are tangent to the plane $SW_{23}$. These points will be called tangent points and can be found from the following equation:

\begin{equation}
	(2,0,-1)
	\begin{pmatrix}
		\frac{a}{3} + \frac{2c}{3}&  b& \frac{2c}{3} - \frac{2a}{3}\\
		-\frac{b}{3}&  a&           \frac{2b}{3}\\
		\frac{c}{3} - \frac{a}{3}& -b&     \frac{2a}{3} + \frac{c}{3}\\
	\end{pmatrix}
	\begin{pmatrix}
		\frac{x_3}{2}-x_{1eq_i}\\x_2\\x_3
	\end{pmatrix}=
	-(a+c)x_{1eq_i}+\frac{3(c-a)}{2}x_3+3bx_2=0,
\end{equation}
it follows that:
\begin{equation}
	x_2=\frac{(a+c)}{3b}x_{1eq_i}-\frac{(c-a)}{2b} x_3, \text{ with } i=1,2.
\end{equation}
For the vector field of $P_2$, $x_{1eq_2}=-\gamma+\alpha$ then:
\begin{equation}\label{eq:tanP2}
	x_2=\frac{(a+c)}{3b}(-\gamma+\alpha)-\frac{(c-a)}{2b} x_3,
\end{equation}
while for the vector field of $P_3$, $x_{1eq_3}=\gamma-\alpha$:
\begin{equation}\label{eq:tanP3}
	x_2=\frac{(a+c)}{3b}(\gamma-\alpha)-\frac{(c-a)}{2b} x_3,
\end{equation}
According to \eqref{eq:surfaces4atoms2}, if $x_3> 0$ then $SW_{23}$ belongs to $P_2$, and the tangent points to consider in $SW_{23}$ for $x_3> 0$ are given by \eqref{eq:tanP2}. And if $x_3\leq 0$, then $SW_{23}$ belongs to $P_3$ and the tangent points are given by \eqref{eq:tanP3}. An illustration of the tangent points in $SW_{23}$ is shown in Figure~\ref{fig:diagrama3d}, where the points for $P_2$ are indicated by a dotted line, while for $P_3$ are drawn as a continuous line.\\

As a starting point to propose the region $R_1$ defined by four points $p_1,\ldots,p_4$, consider the point in $cl(W^u_{{\bf x}_{eq_2}})\cap SW_{23}$ given by \eqref{eq:interseccionp2a} that fulfills \eqref{eq:tanP2}:

\begin{equation}
	pt_1=\begin{pmatrix}
		-\displaystyle\frac{\gamma-\alpha}{3}\\
		\displaystyle-\frac{2a(\gamma-\alpha)}{3b}\\
		\displaystyle-\frac{2(\gamma-\alpha)}{3}
	\end{pmatrix},
\end{equation}
and in ${\bf z}^{(2)}$ coordinates 
\begin{equation}
	pt_1z_2=\begin{pmatrix}
		0\\
		\displaystyle\frac{2a(\gamma-\alpha)}{3b}\\
		\displaystyle-\frac{2(\gamma-\alpha)}{3}
	\end{pmatrix}.
\end{equation}

If we evaluate the trajectory with initial condition in  $x_0=pt_1z_2$, under the vector field of $P_2$ and ignoring the effect of the vector field of $P_1$ and $P_3$, reaches the point $pt_2z_2 \in SW_{23}$. The flow $\varphi$ can scape from $P_2$ to $P_3$ through the segment $\overline{pt_1z_2\;pt_2z_2}$.
Thus, trajectories with initial condition close to ${\bf z}_{eq_2}$ and not in the stable manifold probably cross $SW_{23}$ close to the segment $\overline{pt_1z_2\;pt_2z_2}$ then $R_1$ should include this segment. However, when the vector field of all atoms is considered, trajectories with initial conditions close to $pt_1z_2$ could reach $SW_{23}$ in points whose second component in $z^{(2)}$ coordinates are further from $0$ than the second component in $z^{(2)}$ coordinates of $pt_2z_2$.  This allows us to propose the region $R_1$ based on a larger segment $\overline{pa_1z_2\;pa_2z_2}$ such that $\overline{pt_1z_2\;pt_2z_2}\subset\overline{pa_1z_2\;pa_2z_2}$. Consider the initial condition $paz_2$  given in ${\bf z}^{(2)}$ coordinates by: 
\begin{equation}
	paz_2=\begin{pmatrix}
		0\\0\\\displaystyle-\frac{2(\gamma-\alpha)}{3}\\
	\end{pmatrix},
\end{equation}
then, the radius with respect to ${\bf z}^*_{eq_2}$ would be $\displaystyle\frac{2(\gamma-\alpha)}{3}$. Remember that only the vector field of $P_2$ is considered and the trajectory rotates around the axis $z^{(2)}_1$. Analyzing an impossible case when the trajectory with initial condition in $paz_2$ reaches $SW_{23}$ with an increment of radius that corresponds to $t=2\pi/b$ ($360^o$) then the $z^{(2)}_2$ component of this point is further form $0$ than the $z^{(2)}_2$ component of $pt_2z_2$. Then, we could take $pa_2z_2=pt_1z_2$ and find the $z^{(2)}_2$ component of $pa_1z_2$ from:

\begin{equation}
	\sqrt{\left(\frac{2(\gamma-\alpha)}{3}e^{a\frac{2\pi}{b}}\right)^2-\left(\frac{2(\gamma-\alpha)}{3}\right)^2}.
\end{equation} 

Trying to find simplicity in the argumentation to obtain the $R_1$ region, we consider the following Assumption:
\begin{ass}\label{as:parameters}
	The parameter values fulfill the following relations: $\displaystyle\frac{b}{a}\geq 25$, $2\geq \left\lvert \displaystyle\frac{c}{b}\right\lvert\geq \displaystyle\frac{7}{5}$ and $\displaystyle\frac{\gamma}{\alpha}\geq 10$.
\end{ass}
then:
\begin{equation}
	\sqrt{\left(\frac{2(\gamma-\alpha)}{3}e^{a\frac{2\pi}{b}}\right)^2-\left(\frac{2(\gamma-\alpha)}{3}\right)^2} \leq
	\gamma\sqrt{\left(\frac{2}{3}e^{a\frac{2\pi}{b}}\right)^2-\left(\frac{2}{3}\right)^2}\leq 0.5388\gamma < \frac{3\gamma}{5}.
\end{equation} 

Remember that $z_2^{(2)}=-x_2$, then, the points  $pa_1$ and $pa_2$ are given by 
\begin{equation}
	pa_1=\begin{pmatrix}
		\displaystyle -\frac{\gamma-\alpha}{3}\\
		\displaystyle\frac{3\gamma}{5}\\
		\displaystyle-\frac{2(\gamma-\alpha)}{3}\\
	\end{pmatrix},\
	pa_2=\begin{pmatrix}
		\displaystyle -\frac{\gamma-\alpha}{3}\\
		\displaystyle -\frac{2a(\gamma-\alpha)}{3b}\\
		\displaystyle-\frac{2(\gamma-\alpha)}{3}\\
	\end{pmatrix},
\end{equation}
where $-\displaystyle\frac{2a(\gamma-\alpha)}{3b}$ is the tangent coordinate given by \eqref{eq:tanP2} for $x_3=-\displaystyle\frac{2(\gamma-\alpha)}{3}$.\\

So, let us propose the region $R_1$ delimited by the following four points:

\begin{equation}
	\begin{array}{l}
		p_1=pa_1+\left(\displaystyle\frac{\gamma}{10},0,\displaystyle\frac{\gamma}{5}\right)^T,\\
		p_2=pa_1-\left(\displaystyle\frac{\gamma}{10},0,\displaystyle\frac{\gamma}{5}\right)^T,\\
		p_3=pa_2+\left(-\displaystyle\frac{\gamma}{10},-\displaystyle\frac{(c-a)(-\frac{\gamma}{5})}{2b},-\frac{\gamma}{5}\right)^T,\\
		p_4=pa_2+\left(\displaystyle\frac{\gamma}{10},-\displaystyle\frac{(c-a)(\frac{\gamma}{5})}{2b},\frac{\gamma}{5}\right)^T.\\
	\end{array}
\end{equation}

\begin{figure}
	\centering
	\includegraphics[width=0.5\textwidth]{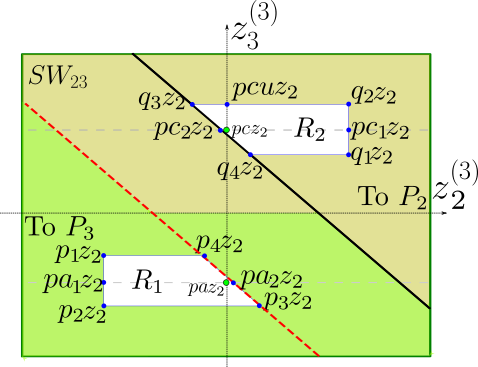}
	\caption{Regions $R_1$ and $R_2$ on the projection $z_2^{(2)}-z_3^{(2)}$.}
	\label{fig:regiones}
\end{figure}

Because $R_2$ is symmetric to $R_1$ with respect to origen, then the symmetric region $R_2$ is delimited by the points:
\begin{equation}
	\begin{array}{l}
		q_1=-p_1,\\
		q_2=-p_2,\\
		q_3=-p_3,\\
		q_4=-p_4.
	\end{array}
\end{equation}
These regions $R_1$ and $R_2$ have been proposed taking into consideration that points $p_d$ and $p_e$ shown in the Figure~\ref{fig:diagramaz2} are part of the regions.
In the Figure~\ref{fig:regiones}, $R_1$ and $R_2$ are shown in ${\bf z}^{(2)}$ coordinates. The points in ${\bf z}^{(2)}$ coordinate system have the suffix $z_2$, for instance $p_2$ in ${\bf z}^{(2)}$ coordinates is $p_2z_2$.\\

The next step is to validate the proposed region $R_1$, the next assumption is made:

\begin{ass}\label{as:rotation}
	Trajectories in $P_1\cup P_2$ rotate around stable manifolds of ${\bf x}_{eq_1}$ and ${\bf x}_{eq_2}$.  The increment in radius for a rotation of $270^o$ is approximately the same if the rotation is performed only around the stable manifold of ${\bf x}_{eq_2}$, {\it i.e.}, if only the vector field of $P_2$ is considered.
\end{ass}

\begin{figure}
	\centering
	\subfloat[]{\includegraphics[width=0.25\columnwidth]{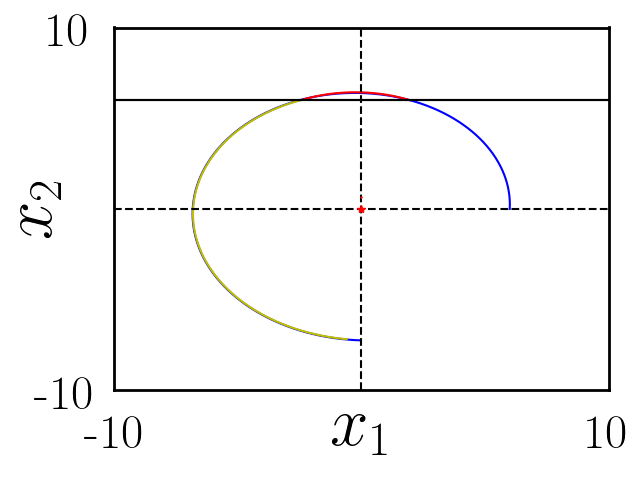}\label{fig:rotation1}}		
	\subfloat[]{\includegraphics[width=0.25\columnwidth]{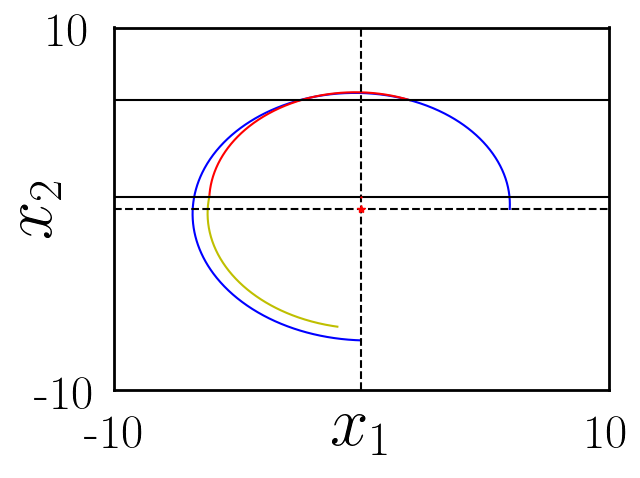}\label{fig:rotation2}}
	\subfloat[]{\includegraphics[width=0.25\columnwidth]{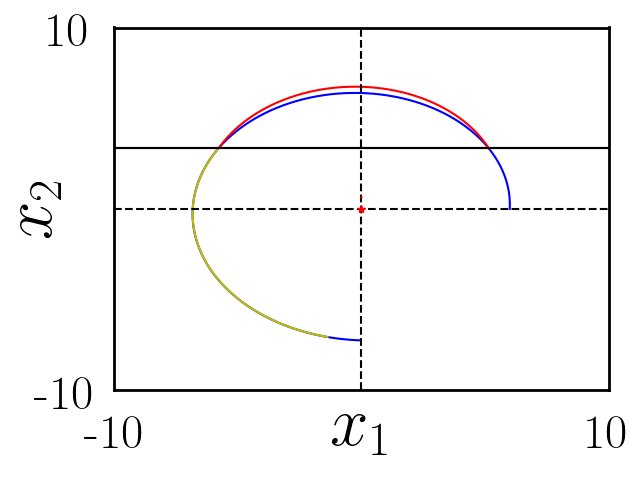}\label{fig:rotation3}}\\
	\subfloat[]{\includegraphics[width=0.25\columnwidth]{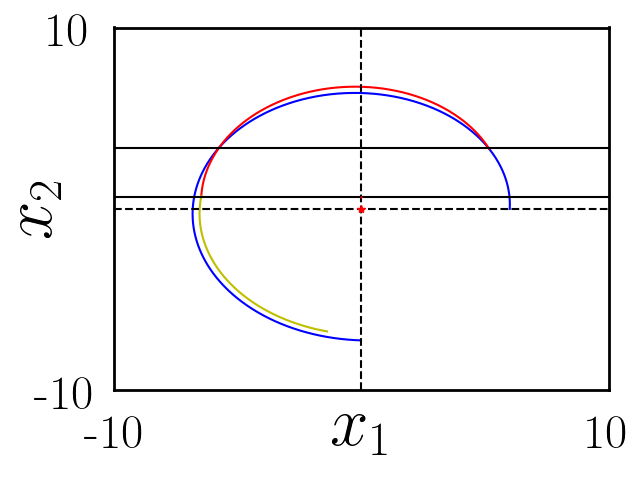}\label{fig:rotation4}}
	\subfloat[]{\includegraphics[width=0.25\columnwidth]{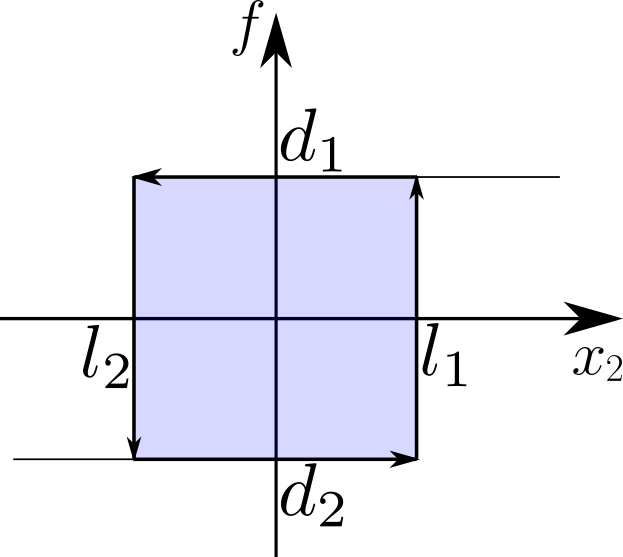}\label{fig:histe}}
	\caption{\label{fig:rotation} Simulation of the system~\eqref{eq:systemahiste} for $t\in[0,\frac{3\pi}{2b})$, $d_1=z^{(2)}_{3eq_1}=\frac{2\alpha}{3}$ and $d_2=z^{(2)}_{3eq_2}=0$, the trajectory for $k=0$, {\it i.e} rotating from the origin is shown in blue. The trajectory in red and yellow is shown for $k=1$ and the parameters (a)$l_1=l_2=\frac{2(\gamma-\alpha)}{3}$, (b)$l_1=\frac{2(\gamma-\alpha)}{3}$, $l_2=\frac{2\alpha}{3}$, (c)$l_1=\frac{\gamma}{3},l_2=\frac{\gamma}{3}$ and (d)$l_1=\frac{\gamma}{3}$, $l_2=\frac{2\alpha}{3}$. In (e) the hysteresis function $f(x_2)$.}
\end{figure}

The assumption \ref{as:rotation} can be easily verified for a set of parameter values, for example,  consider the following system
\begin{equation}\label{eq:systemahiste}
	\dot{x}=\begin{pmatrix}
		a &-b\\b &a
	\end{pmatrix}\begin{pmatrix}
	x_1\\x_2
\end{pmatrix}-k\begin{pmatrix}
a &-b\\b&a
\end{pmatrix}\begin{pmatrix}
0\\f(x_2)
\end{pmatrix},
\end{equation}
where the parameters $a=0.2$, $b=5$, $c=-7$, $\gamma=10$ and $\alpha=1$, and $f(x_2)$ is an hysteresis function as shown in Figure~\ref{fig:histe}

Note that the system~\eqref{eq:systemahiste} is similar to the original system when is projected onto the plane $z_2^{(2)}-z_3^{(2)}$, the commutation at the points $(0,l_1)^T$ and $(0,l_2)^T$ resemble the location of the surface $SW_{23}$ and thus the projection of any trajectory in $P_1\cup P_2$ can be represented for a small duration by the system $\eqref{eq:systemahiste}$.\\

Consider the duration $t\in[0,\frac{3\pi}{2b}]$ and the parameters $d_1=z^{(2)}_{3eq_1}=\frac{2\alpha}{3}$, $d_2=z^{(2)}_{3eq_2}=0$.
In the Figures~\ref{fig:rotation1} - \ref{fig:rotation4},  the trajectory in blue for $k=0$ is shown, {\it i.e}, rotating from the origin.
The trajectory in red and yellow is shown for $k=1$ and the parameters (a)$l_1=l_2=\frac{2(\gamma-\alpha)}{3}$, (b)$l_1=\frac{2(\gamma-\alpha)}{3}$, $l_2=\frac{2\alpha}{3}$, (c)$l_1=\frac{\gamma}{3},l_2=\frac{\gamma}{3}$ and (d)$l_1=\frac{\gamma}{3}$, $l_2=\frac{2\alpha}{3}$.\\

Thus, it can be seen from the simulations of the system~\eqref{eq:systemahiste} for different commutation settings that the assumption \ref{as:rotation} is valid. Moreover, as $\frac{\gamma}{\alpha}$ increases a better approximation is expected.\\

The next step is to verify the validity of the region $R_1$, for this purpose let us define the set $R_1b$ as follows:
\begin{equation}
	R_1b= \left\{ {\bf z}^{(2)}\in\mathbb{R}^3: z_1^{(2)}\in \left[-\frac{\gamma}{5},\frac{\gamma}{5}\right ] \right\}.
\end{equation}
First let us verify that the points in $R_2$ go to $R_1b$.

The evaluation of the vector field in $pc_1z_2$ tell us that the spin is counterclockwise in ${\bf z}^{(2)}$ coordinates. It is not hard to see that the points below the segment $\overline{pcz_2\;pc_1z_2}$ produce trajectories that can perform a turn of $180^o$ around the $z_1^{(2)}$ axis without reaching $SW_{23}$ again. The time that correspond to $360^o$ under the Assumption~\ref{as:rotation} is $T=2\pi/b$. The point $pc_1z_2$ is given by
\begin{equation}
	pc_1z_2=\begin{pmatrix}
		\displaystyle\frac{4(\gamma-\alpha)}{3}\\
		\displaystyle\frac{3\gamma}{5}\\
		0
	\end{pmatrix}.
\end{equation}
Consider the trajectory with initial condition in $pc_1z_2$, and an evolution time that corresponds to a turn of $180^o$ around the $z_1^{(2)}$ axis. After this time the first component of the state vector can be found from:
\begin{equation}
	z_1^{(2)}=\left(\frac{4(\gamma-\alpha)}{3}\right)e^{c\frac{\pi}{b}}.
\end{equation}
If $z_1^{(2)}\leq\gamma/5$ means that the trajectory with initial condition in $pc_1z_2$ reaches $R_1b$. Consider the Assumption~\ref{as:parameters} for a big value of $z_1^{(2)}$  (when $\gamma$ is too big):
\begin{equation}
	\left(\frac{4(\gamma-\alpha)}{3}\right)e^{c\frac{\pi}{b}}\leq
	\frac{4\gamma}{3} e^{c\frac{\pi}{b}}\leq 0.0164\gamma < \frac{\gamma}{5} .
\end{equation}
Thus the set $\{ {\bf z}^{(2)}\in R_2:z_2^{(2)}\geq 0,z_3^{(2)}\leq 0 \}$ reach the set  $R_1b$.\\

Now consider the point $q_2z_2$ given by:
\begin{equation}
	q_2z_2=\begin{pmatrix}
		\displaystyle \frac{23\gamma}{15}-\frac{4\alpha}{3}\\
		\displaystyle\frac{3\gamma}{5}\\
		\displaystyle\frac{\gamma}{10}
	\end{pmatrix}.
\end{equation}
The angle produced by the radius from the point $q_2z_2$ to the $z_1{(2)}$ axis and the plane $z_1^{(2)}-z_2^{(2)}$ is given by:

\begin{equation}
	\left(\frac{360}{2\pi}\right)\arctan\left(\frac{\frac{\gamma}{10}}{\frac{3\gamma}{5}}\right)=\left(\frac{360}{2\pi}\right)\arctan\left(\frac{1}{6}\right)=9.4623^o.
\end{equation}

Let us consider that the trajectory with initial condition in $q_2z_2$ evolves for a duration time that corresponds to $180^o-9.4623^o=170.5377^o$. The first component of the state vector after this duration is given by:
\begin{equation}
	z_1^{(2)}=\left(\frac{23\gamma}{15}-\frac{4\alpha}{3}\right)e^{c\left(\frac{170.5377\pi}{180b}\right)}.
\end{equation}
As before, if $z_1^{(2)}\leq\gamma/5$ means that the trajectory with initial condition in $q_2z_2$ reaches $R_1b$. Consider again the Assumption~\ref{as:parameters} for a big value of $z_1^{(2)}$:
\begin{equation}
	\left(\frac{23\gamma}{15}-\frac{4\alpha}{3}\right)e^{c\left(\frac{170.5377\pi}{180b}\right)}\leq \left(\frac{23\gamma}{15}\right)e^{c\left(\frac{170.5377\pi}{180b}\right)}\leq
	0.0238\gamma < \frac{\gamma}{5}.
\end{equation}
Thus, $q_2z_2$ reaches the region $R_1b$. Moreover, since the points in the segment $\overline{q_2z_2\; pc_1z_2}$ produce radius whose angle with the plane $z_1^{(2)}-z_2^{(2)}$ is between $0$ and $9.4623$ the trajectories starting in this segment also reach the set $R_1b$.\\

Now consider the point $pcuz_2$ given by
\begin{equation}
	pcuz_2=\begin{pmatrix}
		\displaystyle \frac{23\gamma}{15}-\frac{4\alpha}{3}\\
		\displaystyle 0\\
		\displaystyle\frac{\gamma}{10}
	\end{pmatrix}.
\end{equation}
The trajectories with initial condition in the segment $\overline{pcz_2\;pcuz_2}$ can turn $90^o$ without reaching $SW_{23}$. Thus, let us consider that the trajectory with initial condition in $pcuz_2$ evolves for a duration time that corresponds to $90^o$. The first component of the state vector after this duration is given by:
\begin{equation}
	z_1^{(2)}=\left(\frac{23\gamma}{15}-\frac{4\alpha}{3}\right)e^{c\frac{\pi}{2b}}.
\end{equation}

If $z_1^{(2)}\leq\gamma/5$ means that the trajectory with initial condition in $pcuz_2$ reaches $R_1b$. Consider the  Assumption~\ref{as:parameters} for a big value of $z_1^{(2)}$:
\begin{equation}
	\left(\frac{23\gamma}{15}-\frac{4\alpha}{3}\right)e^{c\frac{\pi}{2b}}\leq \left(\frac{23\gamma}{15}\right)e^{c\frac{\pi}{2b}}\leq
	0.17\gamma < \frac{\gamma}{5}.
\end{equation}
Thus the trajectories with initial condition in the set $\left\{ {\bf z}^{(2)}\in R_2:z_2^{(2)}\geq 0,z_3^{(2)}\geq 0 \right\}$ converge to the set  $R_1b$.\\

Now consider the point $q_3z_2$
\begin{equation}
	q_3z_2=\begin{pmatrix}
		\displaystyle \frac{23\gamma}{15}-\frac{4\alpha}{3}\\
		\displaystyle \frac{3c\gamma-23a\gamma}{30b}+\frac{2a\alpha}{3b}\\
		\displaystyle \frac{\gamma}{10}
	\end{pmatrix}.
\end{equation}

To verify that the trajectory starting in $q_3z_2$ is not going to reach $SW_{23}$ when the radius form an angle of $270^o$  with the plane $z_1^{(2)}-z_2^{(2)}$, let us consider the following exaggerated scenario: The radius size corresponds to a duration equivalent to $270^o$ but the $z_1^{(2)}$ component corresponds to a duration equivalent to $90^o$ of oscillation, {\it i.e} when the radius form an angle of $270^o$ with the plane $z_1^{(2)}-z_2^{(2)}$ a smaller radius than the real one is considered, also, a larger value of $z_1^{(2)}$ than the real value is considered.

Then to obtain the radius:

\begin{equation}
	r=e^{a\frac{3\pi}{4b}}\sqrt{\left(\frac{3c\gamma-23a\gamma}{30b}+\frac{2a\alpha}{3b}\right)^2+\left(\frac{\gamma}{10}\right)^2}.
\end{equation}
Under the Assumption~\ref{as:parameters} for the largest radius: 
\begin{equation}
	e^{a\frac{3\pi}{4b}}\sqrt{\left(\frac{3c\gamma-23a\gamma}{30b}+\frac{2a\alpha}{3b}\right)^2+\left(\frac{\gamma}{10}\right)^2}\leq
	e^{a\frac{3\pi}{4b}}\sqrt{\left(\frac{3c\gamma-23a\gamma}{30b}\right)^2+\left(\frac{\gamma}{10}\right)^2},
\end{equation}
\begin{equation}
	e^{a\frac{3\pi}{4b}}\sqrt{\left(\frac{3c\gamma}{30b}-\frac{23a\gamma}{30b}\right)^2+\left(\frac{\gamma}{10}\right)^2}\leq
	e^{a\frac{3\pi}{4b}}\sqrt{\left(-\frac{1\gamma}{5}- \frac{23\gamma}{750}\right)^2+\left(\frac{\gamma}{10}\right)^2},
\end{equation}
\begin{equation}
	e^{a\frac{3\pi}{4b}}\sqrt{\left(-\frac{1\gamma}{5}- \frac{23\gamma}{750}\right)^2+\left(\frac{\gamma}{10}\right)^2}=
	\gamma e^{a\frac{3\pi}{4b}}\sqrt{\left(-\frac{173}{750}\right)^2+(\frac{1}{10})^2},
\end{equation}
\begin{equation}
	\gamma e^{a\frac{3\pi}{4b}}\sqrt{\left(-\frac{173}{750}\right)^2+\left(\frac{1}{10}\right)^2}\leq 0.3036\gamma .
\end{equation}

The value of $z_1^{(2)}$ after a duration equivalent to $90^o$ is given by
\begin{equation}
	z_1^{(2)}=\left(\frac{23\gamma}{15}-\frac{4\alpha}{3}\right)e^{c\frac{\pi}{2b}}.
\end{equation}
Under the Assumption~\ref{as:parameters} for the biggest $z_1^{(2)}$
\begin{equation}
	\left(\frac{23\gamma}{15}-\frac{4\alpha}{3}\right)e^{c\frac{\pi}{2b}}\leq 
	\left(\frac{23\gamma}{15}\right)e^{c\frac{\pi}{2b}},
\end{equation}
\begin{equation} 
	\left(\frac{23\gamma}{15}\right)e^{c\frac{\pi}{2b}}\leq 0.17\gamma < \frac{\gamma}{5}.
\end{equation}
In ${\bf z}^{(2)}$ coordinates, for a specific value of $z_1^{(2)}$  $SW_{23}$ fulfills $z_3^{(2)}=-\displaystyle\frac{2(\gamma-\alpha)}{3}+\frac{z_1^{(2)}}{2}$.
At this angle of $270^o$ the radius is  $r=-z_3^{(2)}$. Then, if the values found for this scenario fulfill the following inequality, it can be concluded that the trajectory with initial condition in $q_3z_2$ does not reach $SW_{23}$ after a duration that corresponds to an oscillation of $90^0$:
\begin{equation}
	-0.3036\gamma>-\frac{2(\gamma-\alpha)}{3}+0.085\gamma ,
\end{equation}
under the Assumption~\ref{as:parameters} for the worst case
\begin{equation}
	-0.3036\gamma>-\frac{18\gamma}{30}+0.085\gamma=-0.515\gamma .
\end{equation}

Then, the trajectory remains for the duration that corresponds to $90^o$. Also, the trajectory does not reach $SW_{23}$ when the radius is at an angle of $270^o$ with the plane  $z_1^{(2)}-z_2^{(2)}$ even when the radius growth is exaggerated. It can be concluded that the trajectory with initial condition in $q_3z_2$ could reach $SW_{23}$ until the second time it approaches $SW_{23}$ and reaches $R_1b$.

Since $q_3z_2$ is the point in the set  $\left\{ {\bf z}^{(2)}\in R_2:z_2^{(2)}\leq 0, \right\} $ that produces the largest radius of that set, the trajectories with initial condition in this set also reach $R_1b$. Then the trajectories starting at $R_2$ reach $R_1b$.

Now to verify that the trajectories that start in $R_2$ reach $R_1\subset R_1b$ is enough to verify the trajectories starting in the segment $\overline{q_1z_2q_2z_2}$, since these produce the largest radius in $R_1b$.\\

Consider the points $q_1z_2$ and $q_2z_2$
\begin{equation}
	q_1z_2=\begin{pmatrix}
		\displaystyle \frac{17\gamma}{15}-\frac{4\alpha}{3}\\
		\displaystyle \frac{3\gamma}{5}\\
		\displaystyle -\frac{\gamma}{10}
	\end{pmatrix},\ 
	q_2z_2=\begin{pmatrix}
		\displaystyle \frac{23\gamma}{15}-\frac{4\alpha}{3}\\
		\displaystyle \frac{3\gamma}{5}\\
		\displaystyle \frac{\gamma}{10}
	\end{pmatrix},
\end{equation}

both points produce the same radius with a different angle, however, more oscillation time before reaching $SW_{23}$ is expected from $q_1z_2$. Thus, consider the trajectory with initial condition in $q_1z_2$ and the evolution time that corresponds to $270^o+9.4623=279.4623$, which is an exaggerated angle since $SW_{23}$ is reached before that.

\begin{equation}
	r=e^{a\frac{279.4623\pi}{180b}}\sqrt{\left(\frac{3\gamma}{5}\right)^2+\left(\frac{\gamma}{10}\right)^2}=\gamma e^{a\frac{279.4623\pi}{180b}}\sqrt{\frac{37}{100}}.
\end{equation}
Under the Assumption~\ref{as:parameters}
\begin{equation}
	\gamma\sqrt{\frac{37}{100}}e^{a\frac{279.4623\pi}{180b}}\leq
	0.7393\gamma.
\end{equation}

\begin{figure}[!ht]
	\centering
	\subfloat[]{\includegraphics[width=0.45\columnwidth]{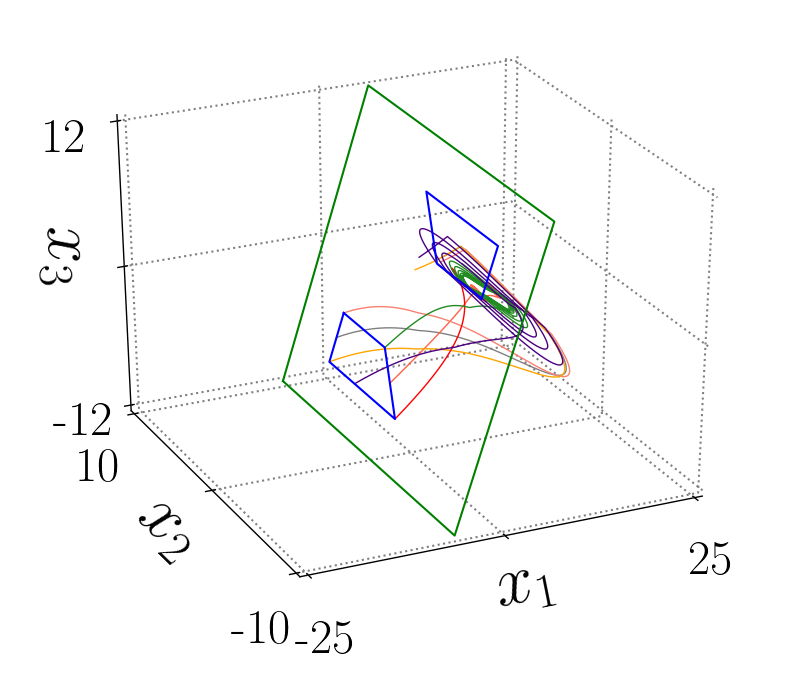}\label{fig:verificacion1}}
	\subfloat[]{\includegraphics[width=0.45\columnwidth]{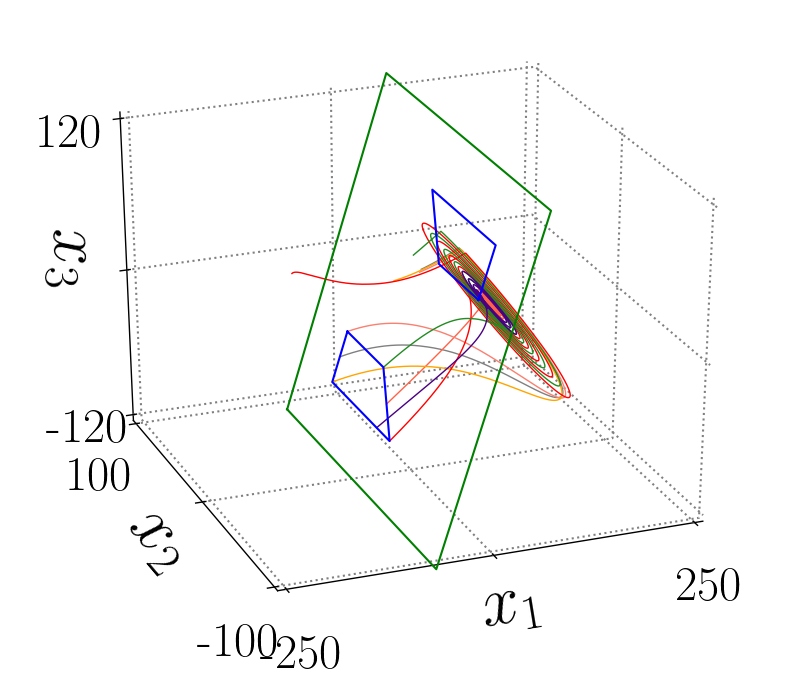}\label{fig:verificacion2}}
	\caption{\label{fig:verificacion}Seven trajectories of the system given by \eqref{eq:affine}, \eqref{eq:Amatrix}, \eqref{eq1:vectorB}, \eqref{eq:funcion4atomos} and \eqref{eq:surfaces4atoms} starting in $R_1$ with $a=0.2$, $b=5$, $c=-7$, $\alpha=1$ and different values of $\gamma$: (a) $\gamma=10$ and (b) $\gamma=100$. }
\end{figure}

Consider the points $p_1z_2$, $pa_1z_2$ and $p_2z_2$
\begin{equation}
	p_1z_2=\begin{pmatrix}
		\displaystyle \frac{\gamma}{5}\\
		\displaystyle -\frac{3\gamma}{5}\\
		\displaystyle \frac{2\alpha}{3}-\frac{17\gamma}{30}
	\end{pmatrix},\
	pa_1z_2=\begin{pmatrix}
		0\\
		\displaystyle -\frac{3\gamma}{5}\\
		\displaystyle -\frac{2(\gamma-\alpha)}{3}
	\end{pmatrix},\ 
	p_2z_2=\begin{pmatrix}
		\displaystyle -\frac{\gamma}{5}\\
		\displaystyle -\frac{3\gamma}{5}\\
		\displaystyle \frac{2\alpha}{3}-\frac{23\gamma}{30}
	\end{pmatrix}.
\end{equation}
The minimum radius in the segment $\overline{pa_1z_2,pa_2z_2}$ is
\begin{equation}
	r=\sqrt{\left(\frac{3\gamma}{5}\right)^2+\left(\frac{2\alpha}{3}-\frac{17\gamma}{30}\right)^2}.
\end{equation}
Under the Assumptions~\ref{as:parameters}, radius from the segment $\overline{pa_1z_2,pa_2z_2}$ is
\begin{equation}
	\sqrt{\left(\frac{3\gamma}{5}\right)^2+\left(\frac{2\alpha}{3}-\frac{17\gamma}{30}\right)^2} \geq 
	\sqrt{\left(\frac{3\gamma}{5}\right)^2+\left(-\frac{\gamma}{2}\right)^2}
\end{equation}
\begin{equation}
	\gamma\sqrt{\left(\frac{3}{5}\right)^2+\left(-\frac{1}{2}\right)^2}=0.781\gamma.
\end{equation}
Since $0.7393\gamma<0.781\gamma$ even when the increment of radius was exaggerated, it can be concluded that the trajectories with initial condition in $R_2$ reach $R_1$ or a self excited attractor. In the same way the trajectories with initial condition in $R_1$ reach $R_2$ or go to a self-excited attractor.\\
To verify the region for the parameters $a=0.2$, $b=5$, $c=-7$ and $\alpha=1$ seven trajectories have been simulated and are shown in the Figure~\ref{fig:verificacion1} for $\gamma=10$ and $\gamma=100$ in the Figure~\ref{fig:verificacion2}.\\ 

Now consider two sets of initial conditions in $SW_{23}$, $I_1$ and $I_2$, such that subsets $N(pa)$ and $N(pc)$ of these sets produce trajectories that end in one of the self excited attractors. These sets are drawn by circles in the Figure~\ref{fig:diagrama3d}.

It is easy to see that if $\gamma$ increases then the regions $R_1$ and $R_2$ grow, but the subsets of initial conditions in $I_1$ and $I_2$ that reach a self-excited attractor without reaching $SW_{23}$ again are reduced.

Let us look at the system in  ${\bf z}^{(2)}$ coordinates, as $\gamma$ grows, $SW_{23}$ and $pcz_2$ is further from the ${\bf z}^*_{eq2}$ and then it takes more time for the trajectories close to $pcz_2$ to travel along the $z^{(2)}_1$ direction to get close to ${\bf z}^*_{eq2}$, however, the expansion along $z^{(2)}_2$ and $z^{(2)}_3$ remains the same, then the subsets of initial conditions that reach the self excited attractors without reaching $SW_{23}$ again shrink in $I_1$ and $I_2$ but never disappear. As $pc$ and $pa$ belong to $R_1$ and $R_2$, respectively, then there will be always an intersection of these regions $R_1$ and $R_2$ with the subsets of initial conditions that reaches the self excited attractors in the regions $I_1$ and $I_2$.

This explains why as $\gamma$ is increased it easy to find initial conditions such that the transitory last long.
Then, to allow the existence of a hidden attractor the intersection of regions $R_1$ and $R_2$ with those sets given by $I_1$ and $I_2$ must be empty, {\it i.e.}, $N(pa)\cap R_1=\emptyset$ and $N(pc)\cap R_2=\emptyset$ .

%%%%%%%%%%%%%
%%Section 5%%
%%%%%%%%%%%%%
\section{Emergence of hidden attractors}\label{Sec_HA}

A way to produce $N(pa)\cap R_1=\emptyset$ and $N(pc)\cap R_2=\emptyset$ and allow the existence of a hidden attractor is by modifying the commutation surface $SW_{23}$ between the two self-excited attractors. Consider the following switching planes:

\begin{equation}\label{eq:SWnueva1}
	\begin{array}{l}
		SW_{12}=cl(P_1)\cup cl(P_2)=\{{\bf x}\in\mathbb{R}^3:2x_1-x_3=-2\gamma,x_1<0\},\\
		SW_{23}=cl(P_2)\cup cl(P_3)=\{{\bf x}\in\mathbb{R}^3:x_1=0\},\\
		SW_{34}=cl(P_3)\cup cl(P_4)=\{{\bf x}\in\mathbb{R}^3:2x_1-x_3=2\gamma,x_1>0\},
	\end{array}
\end{equation}

\begin{figure}[!ht]
	%\centering
	\subfloat[]{\includegraphics[width=0.4\columnwidth]{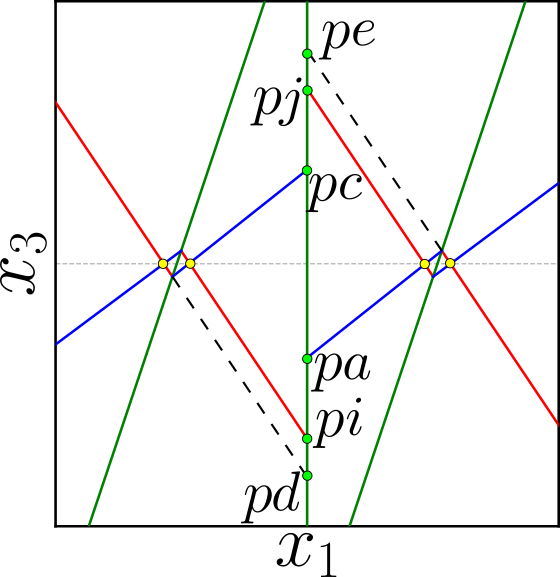}\label{fig:diagramab}}		 \hspace{1cm}
	\subfloat[]{\includegraphics[width=0.45\columnwidth]{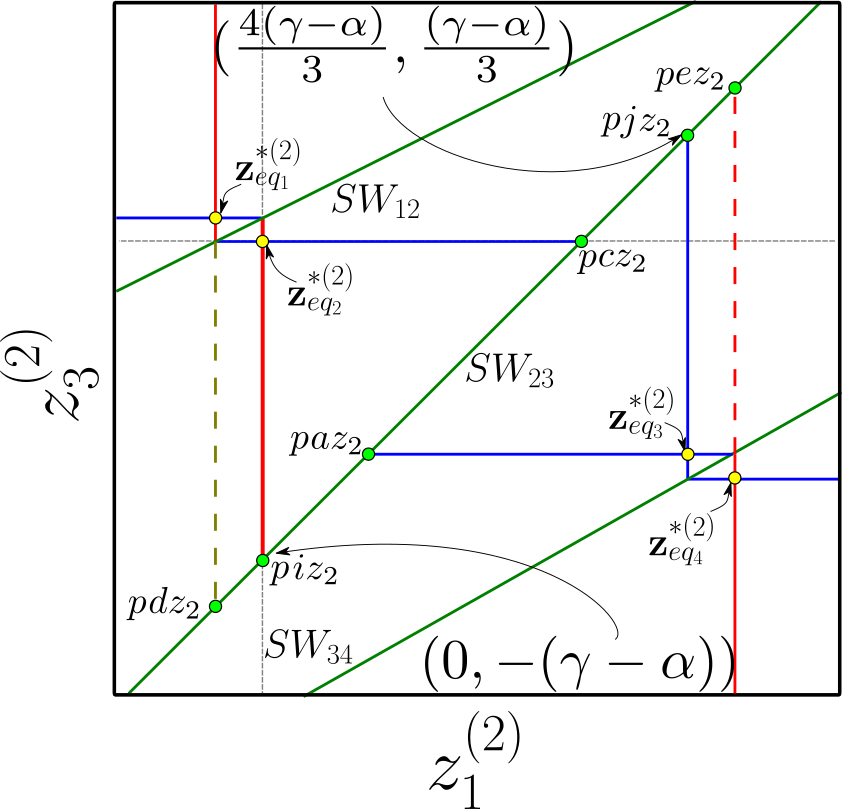}\label{fig:diagramabz2}}
	\hfill
	\caption{\label{fig:diagramasb} Projection of the manifolds on (a)$x_1-x_3$ and (b)$z^{(2)}_1-z^{(2)}_3$. The stable and unstable manifolds are marked with blue and red solid lines, respectively, the switching surfaces with green lines. }
\end{figure}

Note that the switching surface $SW_{23}$ has a new location while the switching surfaces $SW_{12}$ and $SW_{34}$ keep their original locations. This new arrangement keeps the existence of the two heteroclinic loops, and thus the two self-excited attractors. The new projections of the system in ${\bf x}$ and ${\bf z}^{(2)}$ coordinates are shown in the Figure~\ref{fig:diagramasb}. \\

To study the emergence of a hidden attractor, the same procedure of the previous section is followed.

Let us find the points in $SW_{23}$ where the vector field of $P_2$ and $P_3$ are tangent to the plane $SW_{23}$. These points can be found from the following equation:

\begin{equation}
	(1,0,0)
	\begin{pmatrix}
		\frac{a}{3} + \frac{2c}{3}&  b& \frac{2c}{3} - \frac{2a}{3}\\
		-\frac{b}{3}&  a&           \frac{2b}{3}\\
		\frac{c}{3} - \frac{a}{3}& -b&     \frac{2a}{3} + \frac{c}{3}\\
	\end{pmatrix}
	\begin{pmatrix}
		-x_{eq_i}\\x_2\\x_3
	\end{pmatrix}=
	-x_{eq_i}\frac{(a+2c)}{3}+bx_2+\frac{2c-2a}{3}x_3=0, 
\end{equation}

\begin{equation}
	x_2=x_{eq_i}\frac{(a+2c)}{3b}-\frac{2c-2a}{3b}x_3, \text{ with } i=1,2,
\end{equation} 

then for the vector field of $P_2$ we have the expression:
\begin{equation}\label{eq:tangenteP2b}
	x_2=\frac{-(\gamma-\alpha)(a+2c)}{3b}-\frac{2c-2a}{3b}x_3.
\end{equation}

For the vector field of $P_3$ the expression is:
\begin{equation}\label{eq:tangenteP3b}
	x_2=\frac{(\gamma-\alpha)(a+2c)}{3b}-\frac{2c-2a}{3b}x_3.
\end{equation}

Consider the point in $cl(W^u_{{\bf x}_{eq_2}}\cap SW_{23})$ that fulfills \eqref{eq:tangenteP2b}:

\begin{equation}
	pt_1=\begin{pmatrix}
		\displaystyle 0\\
		\displaystyle -\frac{a(\gamma-\alpha)}{b}\\
		\displaystyle -(\gamma-\alpha)
	\end{pmatrix},
\end{equation}

and in ${\bf z}^{(2)}$ coordinates

\begin{equation}
	pt_1z_2=\begin{pmatrix}
		\displaystyle 0\\
		\displaystyle \frac{a(\gamma-\alpha)}{b}\\
		\displaystyle -(\gamma-\alpha)
	\end{pmatrix}.
\end{equation}

If we evaluate the trajectory with initial condition in $x_0=pt_1z_2$, under the vector field of $P_2$ ignoring the effect of the vector field of $P_1$ and $P_3$, reaches the point $pt_2z_2\in SW_{23}$. The flow $\varphi$ can scape from $P_2$ to $P_3$ through the segment $\overline{pt_1z_2pt_2z_2}$. Thus, trajectories with initial condition close to ${\bf z}_{eq_2}$ and not in the stable manifold probably cross $SW_{23}$ close to the segment $\overline{pt_1z_2\;pt_2z_2}$, then, $R_1$ should include this segment. However, when the vector field of all atoms is considered, trajectories with initial conditions close to $pt_1z_2$ could reach $SW_{23}$ in points whose second component in $z^{(2)}$ coordinates are further from $0$ than the second component in $z^{(2)}$ coordinates of $pt_2z_2$.  This allows us to propose the region $R_1$ based on a larger segment $\overline{pi_1z_2\;pi_2z_2}$ such that $\overline{pt_1z_2\;pt_2z_2}\subset\overline{pi_1z_2\;pi_2z_2}$. Consider the initial condition $piz_2$  given in ${\bf z}^{(2)}$ coordinates by: 

\begin{equation}
	piz_2=\begin{pmatrix}
		\displaystyle 0\\
		\displaystyle 0\\
		\displaystyle -(\gamma-\alpha)
	\end{pmatrix},
\end{equation}

then, the radius with respect to ${\bf z}^*_{eq_2}$ would be $(\gamma-\alpha)$. Remember that only the vector field of $P_2$ is considered and the trajectory rotates around the axis $z^{(2)}_1$. Let us analyze an impossible case when the trajectory with initial condition in $piz_2$ reaches $SW_{23}$ with an increment of radius that corresponds to $t=2\pi/b$ ($360^o$). Then, we could take $piz_2=pt_1z_2$ and find the $z_2^{(2)}$ component of $pi_1z_2$ from:

\begin{equation}
	\sqrt{(e^{a\frac{2\pi}{b}}(\gamma-\alpha))^2-((\gamma-\alpha))^2}.
\end{equation} 

Consider the Assumption~\ref{as:parameters}, then:
\begin{equation}
	\sqrt{(e^{a\frac{2\pi}{b}}(\gamma-\alpha))^2-((\gamma-\alpha))^2}\leq
	\gamma\sqrt{(e^{a\frac{2\pi}{b}})^2-1}\leq 0.80815\gamma\approx\frac{4\gamma}{5}.
\end{equation} 

Remember that $z_2^{(2)}=-x_2$, then, the points  $pi_1$ and $pi_2$ are given by 
\begin{equation}
	pi_1=\begin{pmatrix}
		\displaystyle 0\\
		\displaystyle\frac{4\gamma}{5}\\
		\displaystyle-(\gamma-\alpha)\\
	\end{pmatrix},\
	pi_2=\begin{pmatrix}
		\displaystyle 0\\
		\displaystyle -\frac{a(\gamma-\alpha)}{b}\\
		\displaystyle -(\gamma-\alpha)\\
	\end{pmatrix},
\end{equation}
where $-\frac{a(\gamma-\alpha)}{b}$ is the tangent coordinate given by \eqref{eq:tangenteP2b} for $x_3=-(\gamma-\alpha)$.\\

Let us propose a region $R_1$ delimited by the following four points:
\begin{equation}
	\begin{array}{l}
		p_1=pi_1+\left(\displaystyle 0,0,\frac{\gamma}{5} \right)^T,\\
		p_2=pi_1-\left(\displaystyle0,0,\frac{\gamma}{5}\right)^T,\\
		p_3=pi_2+\left(\displaystyle0,-\left(\frac{2c-2a}{3b}\right)\left(-\frac{\gamma}{5}\right),-\frac{\gamma}{5}\right)^T,\\
		p_4=pi_2+\left(\displaystyle0,-\left(\frac{2c-2a}{3b}\right)\left(\frac{\gamma}{5}\right),\frac{\gamma}{5}\right)^T.
	\end{array}
\end{equation}

\begin{figure}
	\centering
	\includegraphics[width=0.5\textwidth]{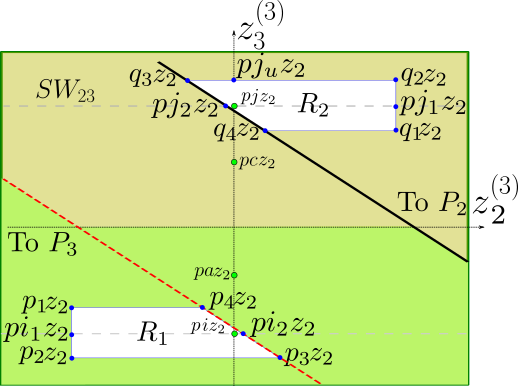}
	\caption{Regions $R_1$ and $R_2$ on the projection $z_2^{(2)}-z_3^{(2)}$.}
	\label{fig:regionesb}
\end{figure}

Then the symmetric region $R_2$ is delimited by the points:
\begin{equation}
	\begin{array}{l}
		q_1=-p_1,\\
		q_2=-p_2,\\
		q_3=-p_3,\\
		q_4=-p_4.
	\end{array}
\end{equation}
The regions $R_1$ and $R_2$ have been proposed taking into consideration that $pd$ and $pe$ are part of the regions and $pa$ and $pc$ are not. In the Figure~\ref{fig:regionesb}, $R_1$ and $R_2$ are shown in ${\bf z}^{(2)}$ coordinates. The points in ${\bf z}^{(2)}$ coordinate system have the suffix $z_2$, for instance $p_2$ in ${\bf z}^{(2)}$ coordinates is $p_2z_2$.\\

The next step is to verify the validity of the region $R_1$, since $p_1z_2-pi$ and $piz_2-p_2z_2$ are given by
\begin{equation}
	p_1z_2-pi=\begin{pmatrix}
		\displaystyle \frac{2\gamma}{15}\\
		\displaystyle -\frac{4\gamma}{5}\\
		\displaystyle \frac{2\gamma}{15}
	\end{pmatrix}, \ 
	pi-p_2z_2=\begin{pmatrix}
		\displaystyle \frac{2\gamma}{15}\\
		\displaystyle \frac{4\gamma}{5}\\
		\displaystyle \frac{2\gamma}{15}
	\end{pmatrix},
\end{equation}
let us define the set $R_1b$ as follows:

\begin{equation}
	R_1b= \left\{ {\bf z}^{(2)}\in\mathbb{R}^3: z_1^{(2)}\in \left[-\frac{2\gamma}{15},\frac{2\gamma}{15}\right] \right\}.
\end{equation}
First let us verify that the points in $R_2$ go to $R_1b$.
The evaluation of the vector field in $pj_1z_2$ tell us that the spin is counterclockwise in ${\bf z}^{(2)}$ coordinates.

Consider the point $q_2z_2$ given by:
\begin{equation}
	q_2z_2=\begin{pmatrix}
		\displaystyle   \frac{22\gamma}{15}-\frac{4\alpha}{3}\\
		\displaystyle\frac{4\gamma}{5}\\
		\displaystyle\frac{7\gamma}{15}-\frac{\alpha}{3}
	\end{pmatrix}.
\end{equation}
The angle produced by the radius from the point $q_2z_2$ to the $z_1{(2)}$ axis and the plane $z_1^{(2)}-z_2^{(2)}$ is given by:

\begin{equation}
	\left(\frac{360}{2\pi}\right)\arctan\left(\frac{\displaystyle\frac{7\gamma}{15}-\frac{\alpha}{3}}{\frac{4\gamma}{5}}\right)=\left(\frac{360}{2\pi}\right)\arctan\left(\frac{7}{12}-\frac{5\alpha}{12\gamma}\right).
\end{equation}

Under the Assumption~\ref{as:parameters} the angle obey the following inequality
\begin{equation}
	\left(\frac{360}{2\pi}\right)\arctan\left(\frac{7}{12}-\frac{5\alpha}{12\gamma}\right)\leq \left(\frac{360}{2\pi}\right)\arctan\left(\frac{7}{12}\right)\leq 30.2564^o.
\end{equation}

Consider then that the trajectory with initial condition in $q_2z_2$ evolves for a duration that corresponds to $180^o-30.2564^o=149.7435^o$. If after this duration the first component of the state vector $z_1^{(2)}\leq2\gamma/15$ means that the trajectory with initial condition in $q_2z_2$ reaches $R_1b$.
After this duration $z_1^{(2)}$ is given by
\begin{equation}
	z_1^{(2)}=\left(\frac{22\gamma}{15}-\frac{4\alpha}{3}\right)e^{c\frac{149.7435\pi}{180b}}.
\end{equation}
Under the Assumption~\ref{as:parameters} for a big value of $z_1^{(2)}$:
\begin{equation}
	\left(\frac{22\gamma}{15}-\frac{4\alpha}{3}\right)e^{c\frac{149.7435\pi}{180b}}\leq \left(\frac{22\gamma}{15}\right)e^{c\frac{149.7435\pi}{180b}}\leq
	0.0378\gamma < \frac{2\gamma}{15}.
\end{equation}
Thus, $q_2z_2$ reaches the region $R_1b$. Moreover since the points in the segment $\overline{q_2z_2q_1z_2}$ produce radius whose angle with the plane $z_1^{(2)}-z_2^{(2)}$ is between $0$ and $30.2564$ the trajectories with initial condition in this segment also reach the set $R_1b$.\\

Now consider the point $q_3z_2$
\begin{equation}
	q_3z_2=\begin{pmatrix}
		\displaystyle \frac{22\gamma}{15}-\frac{4\alpha}{3}\\
		\displaystyle \frac{2c\gamma-17a\gamma}{15b}+\frac{a\alpha}{b}\\
		\displaystyle \frac{7\gamma-5\alpha}{15}
	\end{pmatrix}.
\end{equation}
The angle of the radius at $q_3z_2$ with the plane $z_1^{(2)}-z_2^{(2)}$ is
\begin{equation}
	-\left(\frac{360}{2\pi}\right)\arctan\left(\frac{\frac{2c\gamma-17a\gamma}{15b}+\frac{a\alpha}{b}}{\frac{7\gamma-5\alpha}{15}}\right)=
	-\left(\frac{360}{2\pi}\right)\arctan\left(\frac{2c\gamma-17a\gamma+15a\alpha}{7b\gamma-5b\alpha}\right).
\end{equation}

Under Assumption~\ref{as:parameters} the angle should be less than:
\begin{equation}
	-\left(\frac{360}{2\pi}\right)\arctan\left(\frac{2c\gamma-17a\gamma+15a\alpha}{7b\gamma-5b\alpha}\right)\leq
	-\left(\frac{360}{2\pi}\right)\arctan\left(
	\frac{2c}{6.5b}-\frac{17a}{6.5b}+\frac{1.5a}{6.5b}
	\right).
\end{equation}
\begin{equation}
	-\left(\frac{360}{2\pi}\right)\arctan\left(
	\frac{2c}{6.5b}-\frac{17a}{6.5b}+\frac{1.5a}{6.5b}
	\right) \leq
	-\left(\frac{360}{2\pi}\right)\arctan\left(
	\frac{4}{6.5}-\frac{17}{6.5(25)}+\frac{1.5}{6.5(25)}
	\right)=35.4^o.
\end{equation}

To verify that the trajectory starting in $q_3z_2$ is not going to reach $SW_{23}$ when the radius form an angle of $270^o$  with the plane $z_1^{(2)}-z_2^{(2)}$, let us consider the following exaggerated scenario: The radius size corresponds to a duration equivalent to $200^o$ but the $z_1^{(2)}$ component corresponds to a duration equivalent to $120^o$ of oscillation, {\it i.e} when the radius form an angle of $270^o$ with the plane $z_1^{(2)}-z_2^{(2)}$ a smaller radius than the real one is considered, also, a larger value of $z_1^{(2)}$ than the real value is considered.

Then to obtain the radius:

\begin{equation}
	r=e^{a\frac{20\pi}{18b}}\sqrt{\left(\frac{2c\gamma-17a\gamma}{15b}+\frac{a\alpha}{b}\right)^2+\left(\frac{7\gamma-5\alpha}{15}\right)^2}.
\end{equation}
Under the Assumption~\ref{as:parameters} for the largest radius
\begin{equation}
	e^{a\frac{20\pi}{18b}}\sqrt{\left(\frac{2c\gamma-17a\gamma}{15b}+\frac{a\alpha}{b}\right)^2+\left(\frac{7\gamma-5\alpha}{15}\right)^2}\leq
	e^{a\frac{20\pi}{18b}}\sqrt{\left(\frac{2c\gamma-17a\gamma}{15b}\right)^2+\left(\frac{7\gamma}{15}\right)^2},
\end{equation}
\begin{equation}
	e^{a\frac{20\pi}{18b}}\gamma\sqrt{\left(\frac{2c-17a}{15b}\right)^2+\left(\frac{7}{15}\right)^2}\leq 0.5984\gamma.
\end{equation}

The value of $z_1^{(2)}$ after the duration that corresponds to $120^o$ is:
\begin{equation}
	\left(\frac{22\gamma}{15}-\frac{4\alpha}{3}\right)e^{c\frac{12\pi}{18b}}.
\end{equation}
Under the Assumption~\ref{as:parameters} for the biggest value of $z_1^{(2)}$
\begin{equation}
	\left(\frac{22\gamma}{15}-\frac{4\alpha}{3}\right)e^{c\frac{12\pi}{18b}}\leq 
	\left(\frac{22\gamma}{15}\right)e^{c\frac{12\pi}{18b}},
\end{equation}
\begin{equation} 
	\left(\frac{22\gamma}{15}\right)e^{c\frac{12\pi}{18b}}\leq0.07814\gamma<\frac{2\gamma}{15}.
\end{equation}
The points in $SW_{23}$ fulfills the following equation:
\begin{equation}
	z_3^{(2)}=-(\gamma-\alpha)+z_1^{(2)}.
\end{equation} 
At this angle of $270^o$ the radius is  $r=-z_3^{(2)}$. Then, if the found  values for this scenario fulfill the following inequality, it can be concluded that the trajectory with initial condition in $q_3z_2$ does not reach $SW_{23}$ after a duration that corresponds to an oscillation of $120^0$:
\begin{equation}
	-.5984\gamma>-(\gamma-\alpha)+0.07814\gamma,
\end{equation}
under the Assumption~\ref{as:parameters} and the worst case
\begin{equation}
	-.5984\gamma>-\frac{9\gamma}{10}+0.07814\gamma=-0.82186\gamma.
\end{equation}

Then, the trajectory remains for the duration that corresponds to $120^o$. Furthermore, since $z_1^{(2)}<2\gamma/15$ the trajectory reaches $R_1b$.
Since $q_3z_2$ is the point in the set  $\left\{ {\bf z}^{(2)}\in R_2:z_2^{(2)}\leq 0 \right\} $ that produces the largest radius of that set, the same conclusion applies for the points in this set.\\
For the trajectories with initial condition in the set $\left\{ {\bf z}^{(2)}\in R_2:z_2^{(2)}> 0 \right\} $ the duration is the equivalent to more than $120^o$. 
Thus, the trajectories that start in $R_2$ reach the set $R_1b$.\\

Now to verify that the trajectories that start in $R_2$ reach $R_1\subset R_1b$ is enough to verify the trajectories starting in the segment $\overline{q_1z_2q_2z_2}$, since these produce the largest radius in $R_1b$.\\

Consider the points $q_1z_2$ and $q_2z_2$
\begin{equation}
	q_1z_2=\begin{pmatrix}
		\displaystyle \frac{6\gamma}{5}-\frac{4\alpha}{3}\\
		\displaystyle \frac{4\gamma}{5}\\
		\displaystyle \frac{\gamma}{5}-\frac{\alpha}{3}
	\end{pmatrix},\ 
	q_2z_2=\begin{pmatrix}
		\displaystyle \frac{22\gamma}{15}-\frac{4\alpha}{3}\\
		\displaystyle \frac{4\gamma}{5}\\
		\displaystyle \frac{7\gamma}{15}-\frac{\alpha}{3}
	\end{pmatrix},
\end{equation}

The largest radius is at $q_2z_2$ while the smaller angle is at $q_1z_2$, then let us consider that radius of $q_2z_2$ with the angle of $q_1z_2$ and the end position at $270^o$ with respect to the plane $z_1^{(2)}-z_2^{(2)}$ which is more that the possible rotation.
The angle is given by
\begin{equation}
	\left(\frac{360}{2\pi}\right)\arctan\left(\frac{\frac{\gamma}{5}-\frac{\alpha}{3}}{\frac{4\gamma}{5}}\right)=\left(\frac{360}{2\pi}\right)\arctan\left(\frac{1}{4}- \frac{5\alpha}{12\gamma}\right).
\end{equation}
Under the Assumption~\ref{as:parameters} the smallest angle is
\begin{equation}
	\left(\frac{360}{2\pi}\right)\arctan\left(\frac{1}{4}- \frac{5\alpha}{12\gamma}\right)\geq 
	\left(\frac{360}{2\pi}\right)\arctan\left( \frac{5}{24}\right)=11.7682^o,
\end{equation}
thus for $270-11.7682=258.2317$
\begin{equation}
	\sqrt{\left(\frac{4\gamma}{5}\right)^2+\left(\frac{7\gamma}{15}-\frac{\alpha}{3}\right)^2}e^{a\frac{258.2317\pi}{180b}}.
\end{equation}
Under the Assumption~\ref{as:parameters} for the biggest radius
\begin{equation}
	\sqrt{\left(\frac{4\gamma}{5}\right)^2+\left(\frac{7\gamma}{15}-\frac{\alpha}{3}\right)^2}e^{a\frac{258.2317\pi}{180b}}
	\leq
	\sqrt{\left(\frac{4\gamma}{5}\right)^2+\left(\frac{7\gamma}{15}\right)^2}e^{a\frac{258.2317\pi}{180b}},
\end{equation}
\begin{equation}
	\sqrt{\left(\frac{4\gamma}{5}\right)^2+\left(\frac{7\gamma}{15}\right)^2}e^{a\frac{258.2317\pi}{180b}}\leq
	\gamma\sqrt{\frac{193}{225}}e^{a\frac{258.2317\pi}{180b}}\leq
	1.1091\gamma.
\end{equation}

Consider the points $p_1z_2$, $pa_1z_2$ and $p_2z_2$
\begin{equation}
	p_1z_2=\begin{pmatrix}
		\displaystyle \frac{2\gamma}{15}\\
		\displaystyle -\frac{4\gamma}{5}\\
		\displaystyle \alpha-\frac{13\gamma}{15}
	\end{pmatrix},\
	pa_1z_2=\begin{pmatrix}
		0\\
		\displaystyle -\frac{4\gamma}{5}\\
		\displaystyle \alpha-\gamma
	\end{pmatrix},\ 
	p_2z_2=\begin{pmatrix}
		\displaystyle -\frac{2\gamma}{15}\\
		\displaystyle -\frac{4\gamma}{5}\\
		\displaystyle \alpha-\frac{17\gamma}{15}
	\end{pmatrix}.
\end{equation}
The minimum radius in the segment $\overline{pa_1z_2,pa_2z_2}$ is
\begin{equation}
	\sqrt{\left(\frac{4\gamma}{5}\right)^2+\left(\alpha-\frac{13\gamma}{15}\right)^2}.
\end{equation}
Under the Assumptions~\ref{as:parameters} for the smallest radius:
\begin{equation}
	\sqrt{\left(\frac{4\gamma}{5}\right)^2+\left(\alpha-\frac{13\gamma}{15}\right)^2}\geq
	\sqrt{\left(\frac{4\gamma}{5}\right)^2+\left(\frac{23\gamma}{30}\right)^2}=1.1081\gamma.
\end{equation}

$1.1092\gamma\approx 1.1081\gamma$, even when the angle of rotation was exaggerated, it can be concluded that the trajectories with initial condition in $R_1$ reach $R_2$ or a self excited attractor. In the same way the trajectories with initial condition in $R_2$ reach $R_1$ or go to a self-excited attractor.\\

Let us look at the system in  ${\bf z}^{(2)}$ coordinates, as in the previous section with the previous  switching surfaces, as $\gamma$ grows, $SW_{23}$ and $pcz_2$ is further from the ${\bf z}^*_{eq2}$ and then it take more time for the trajectories close to $pcz_2$ to travel along the $z^{(2)}_1$ direction to get close to ${\bf z}^*_{eq2}$, however, the expansion along $z^{(2)}_2$ and $z^{(2)}_3$ remains the same, then the subsets of initial conditions that reach the self excited attractors without reaching $SW_{23}$ again shrink in $I_1$ and $I_2$ but this time, as opposed to the previous case there exist a value of $\gamma$ such that the intersection disappear.\\

Then for a sufficiently big value of $\gamma$ we have a region $R_1$ such that any trajectory starting there remains crossing $R_1$ for $t>0$. Then, there exists either a periodic orbit, a hidden limit cycle, a hidden chaotic attractor or a combination of the previous, which should go trough $R_1$ and $R_2$.

Also, as small differences in the initial conditions in $R_1$ could eventually produce a big separation of trajectories in $SW_{23}$, sensitivity to initial conditions could also be expected, however the formal proof is out of the scope of this work.
\begin{figure}[!ht]
	\centering
	\subfloat[]{\includegraphics[width=0.45\columnwidth]{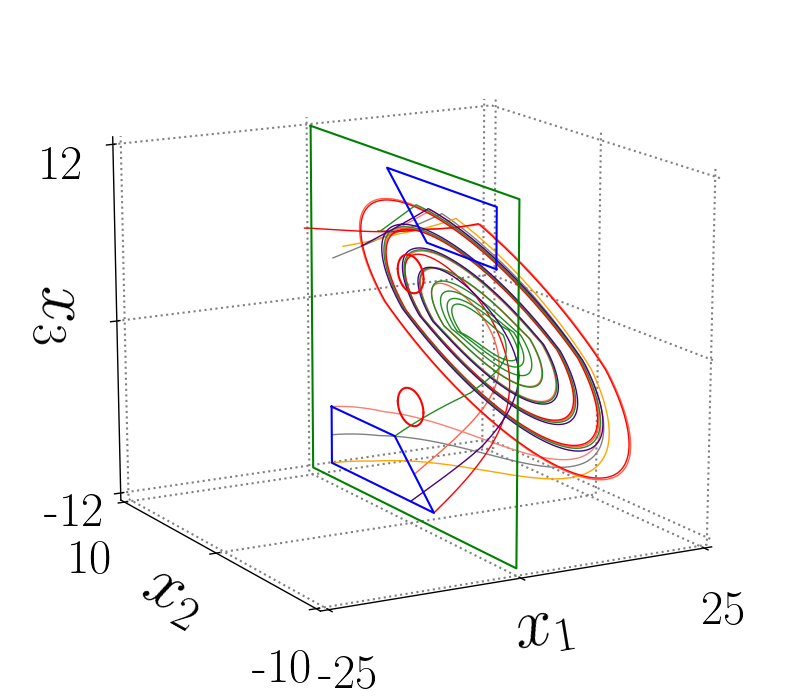}\label{fig:verificacion3}}
	\subfloat[]{\includegraphics[width=0.45\columnwidth]{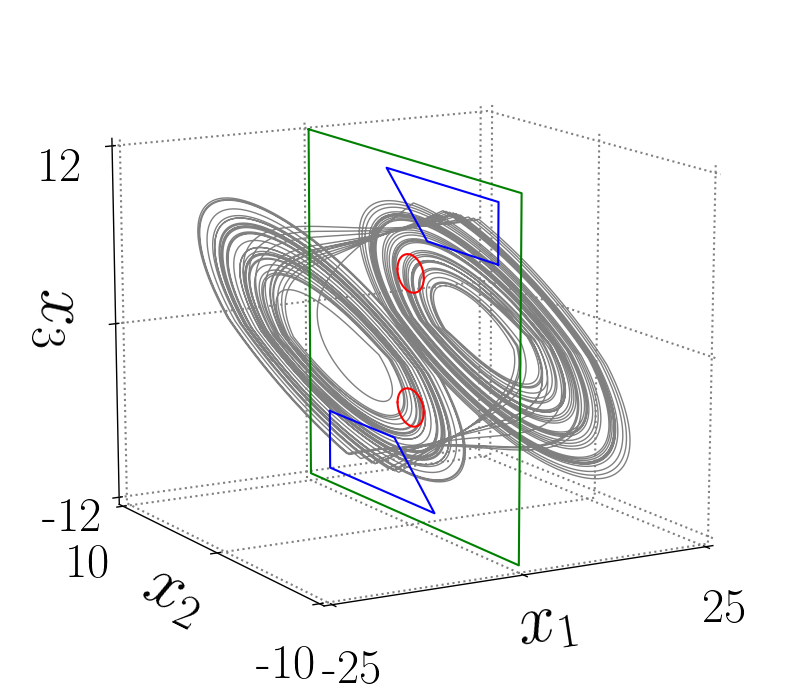}\label{fig:verificacion4}}\\
	\subfloat[]{\includegraphics[width=0.45\columnwidth]{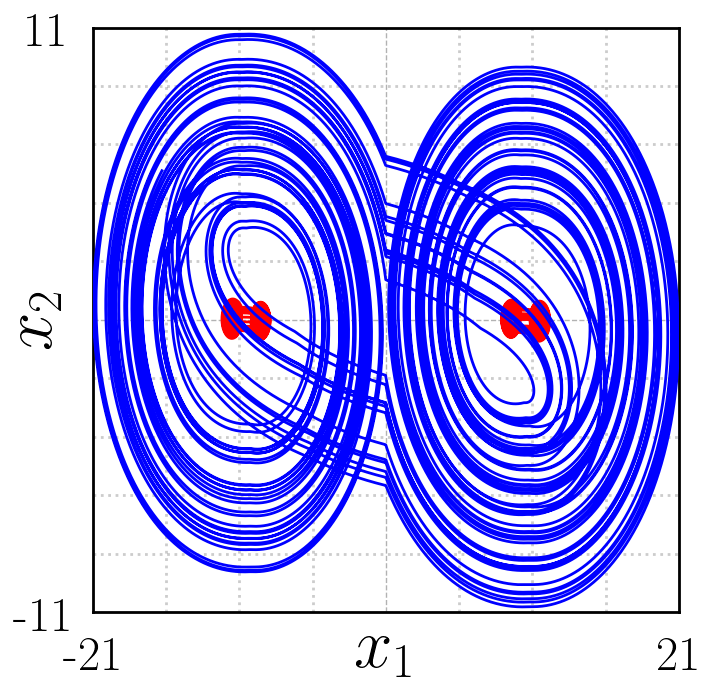}\label{fig:verificacion5}}
	\caption{\label{fig:verificacionb}In (a) seven trajectories of the system given by \eqref{eq:affine}, \eqref{eq:Amatrix}, \eqref{eq1:vectorB}, \eqref{eq:funcion4atomos} and \eqref{eq:SWnueva1} starting in $R_1$ with $a=0.2$, $b=5$, $c=-7$, $\alpha=1$ and $\gamma=10$. In (b) hidden attractor for the same parameters and the initial condition ${\bf x}_0=(0,0,0)^T$ for $t\in[50000,50100]$. In (c) the projection of the self-excited attractors and the hidden attractor onto the plane $x_1-x_2$.} 
\end{figure}

To verify the region for the parameters $a=0.2$, $b=5$, $c=-7$, $\alpha=1$ and $\gamma=10$ seven trajectories have been simulated and are shown in the Figure~\ref{fig:verificacion3}.

Simulation experiments verify the conjecture on the emergence of the hidden attractor.
In the Figure~\ref{fig:verificacion4} it is shown the hidden attractor for the parameters $a=0.2$, $b=5$, $c=-7$, $\alpha=1$ and $\gamma=10$ and the initial condition ${\bf x}_0=(0,0,0)^T$ for $t\in[50000,50100]$.
In the Figure~\ref{fig:verificacion5} it is shown the projection of the hidden attractor and the two self-excited attractors onto the plane $x_1-x_2$ for the same parameters.

\section{Conclusions}\label{Sec:conclusion}
In this work an approach for the generation of multiscroll attractors was studied based on heteroclinic orbits. Particularly, we presented a quad-scroll self-excited attractor which is split in two double-scroll self-excited attractor, so the system bifurcates from monostability to biestability.  The approach is based on the coexistence of double scroll self-excited attractors surrounded the equilibria and presenting heteroclinic orbits. Increasing the distances between the double-scroll self-excited attractors generates a heteroclinic-like orbit between equilibria of two different double-scroll self-excited attractors. It is possible to generate hidden attractors surrounded the self-excited attractors by breaking the heteroclinic-like orbit. The study revealed a relation between the existence of a hidden attractor and the trajectories that resemble heteroclinic orbits at a larger scale which join the self-excited attractors. The findings suggest that new classes of  multistable systems with different number of self-excited and hidden attractors can be designed with a geometric approach.

\section*{Acknowledgments}
R.J.Escalante-Gonz\'alez is thankful to CONACYT (M\'exico) for the scholarships granted. E. Campos-Cant\'on acknowledges CONACYT for the financial support through Project No. A1-S-30433.

\end{document}